\documentclass[12pt,a4paper] {article}
\usepackage{amssymb,amsthm}
\usepackage{graphics,graphicx}
\usepackage[T1]{fontenc}
\usepackage{color}

\def\EE{\mathbb{E}}
\def\ZZ{\mathbb{Z}}
\def\RR{\mathbb{R}}
\def\TT{\mathbb{T}}
\def\TT{\mathbb{T}}
\def\NN{\mathbb{N}}
\def\PP{\mathbb{P}}

\def\exp{{\rm exp}}

\def\Proof {\vskip -2mm {\underline {Proof}}}
\def\fdm {\hfill\break \boxit{2pt} {}}
\def\build#1_#2^#3{\mathrel{\mathop{\kern 0pt#1}\limits_{#2}^{#3}}}

\font\fivegoth=eufm5 \font\sevengoth=eufm7 \font\tengoth=eufm10
\newfam\gothfam
\scriptscriptfont\gothfam=\fivegoth \textfont\gothfam=\tengoth
\scriptfont\gothfam=\sevengoth

\def\boxit#1#2{\setbox1=\hbox{\kern#1{#2}\kern#1}%
\dimen1=\ht1 \advance\dimen1 by #1 \dimen2=\dp1 \advance\dimen2 by#1
\setbox1=\hbox{\vrule height\dimen1 depth\dimen2\box1\vrule}%
\setbox1=\vbox{\hrule\box1\hrule}%
\advance\dimen1 by .4pt \ht1=\dimen1 \advance\dimen2 by .4pt
\dp1=\dimen2 \box1\relax}
\def \fdm {\boxit{3pt} {}}

\newtheorem{hypo}{Hypothesis}[section]
\newtheorem{cond}[hypo]{Condition}

\newtheorem{ppty}[hypo]{Property}

\newtheorem{prop}[hypo]{Proposition}
\newtheorem{thm}[hypo]{Theorem}
\newtheorem{lem}[hypo]{Lemma}

\newtheorem{notas}[hypo]{Notations}

\newtheorem{rem}[hypo]{Remark}
\newtheorem{rems}[hypo]{Remarks}
\newtheorem{coro}[hypo]{Corollary}

\begin{document}
\parindent=0mm
\author{J.-P. Conze, S. Le Borgne, M. Roger
\footnote{IRMAR, UMR CNRS 6625, Universit\'e de Rennes I, Campus
de Beaulieu, 35042 Rennes Cedex, France}}
\title{Central limit theorem for products  \\ of
toral automorphisms}

\date{11 June 2010}

\markright{Central limit theorem for products}


\maketitle

{\abstract Let $(\tau_n)$ be a sequence of toral automorphisms
$\tau_n : x \rightarrow A_n x \hbox{ mod }\ZZ^d$ with $A_n \in {\cal
A}$, where ${\cal A}$ is a finite set of matrices in $SL(d,
\mathbb{Z})$. Under some conditions the method of "multiplicative
systems" of Koml\`os can be used to prove a Central Limit Theorem
for the sums $\sum_{k=1}^n f(\tau_k \circ \tau_{k-1} \cdots \circ
\tau_1 x)$ if $f$ is a H\"older function on $\mathbb{T}^d$. These
conditions hold for $2\times 2$ matrices with positive coefficients.
In dimension $d$ they can be applied when $A_n= A_n(\omega)$, with
independent choices of $A_n(\omega)$ in a finite set of matrices
$\in SL(d, \mathbb{Z})$, in order to prove a "quenched" CLT.}


\section*{Introduction}

Let us consider a sequence of maps obtained by composition of
transformations $(\tau_n)$ acting on a probability space $(X, {\cal
B}, \lambda)$. The iteration of a single measure preserving transformation
corresponds to the classical case of a dynamical system. The case of
several transformations has been also considered by some authors,
and the stochastic behavior of the sums $\sum_{k=1}^n f(\tau_k \circ
\tau_{k-1} \cdots \circ \tau_1 x)$, for a function $f$ on $X$, has
been studied on some examples. For example the notion of stochastic
stability is defined using composition of transformations chosen at
random in the neighborhood of a given one. Bakhtin considered in
\cite{Ba95} non perturbative cases with geometrical assumptions on
the transformations. In the non-invertible case, the example of
sequences of expanding maps of the interval was carried out in
\cite{CoRa07}.

\vskip 3mm Here we consider the example of automorphisms of the
torus. Given a finite set ${\cal A}$ of matrices in
$SL(d,\mathbb{Z})$, to a sequence $(A_{i})_{i\in\mathbb{N}}$ taking
values in ${\cal A}$ corresponds the sequence
$(\tau_{i})_{i\in\mathbb{N}}$ of automorphisms of the torus
$\mathbb{T}^d$ defined by: $\tau_{i} : x\mapsto A_{i}^t x \textrm{
mod } 1$. If the choice in  ${\cal A}$ of the matrices is random, we
write $A_i(\omega)$ and $\tau_{i}(\omega)$.

Let $f:\mathbb{T}^d\rightarrow\mathbb{R}$ be a H\"older function with
integral zero. A question is the existence of the variance and the
central limit theorem for the sums $S_N f=\sum_{k=1}^{N}f(\tau_k...
\tau_1.)$ and the Lebesgue measure $\lambda$ on the torus.

When the matrices are chosen at random and independently, our
problem is strongly related to the properties of a random walk on
$SL(d, \ZZ)$. Among the many works on random walks on groups let us
mention a paper of Furman and Shalom which deals with questions
directly connected to ours. Let $\mu$ be a probability measure on
$SL(d, \mathbb{Z})$. Let $\PP =\mu^{\otimes \NN}$ the product measure on $\Omega:=SL(d, \mathbb{Z})^{\NN}$. In \cite{FuSh99},
if the group generated by the support of $\mu$ has no abelian
subgroup of finite index and acts irreducibly on $\RR^d$, a spectral
gap is proved for the convolution by $\mu$ on $L^2_0$ and a CLT is
deduced for $f(\tau_k(\omega)\ldots \tau_1(\omega)x)$ as a random
variable defined on $(\Omega\times\TT^d, \PP\otimes \lambda)$.
Remark that results of Derriennic and Lin \cite{{DeLi03}} imply the
CLT for $f$ in $L^2_0$ not only for the stationary measure of the
Markov chain, but also for $\lambda$-almost every $x$, with respect
to the measure starting from $x$. This is a quenched CLT, but with a
meaning different from ours: for them $x$ is fixed, for us $\omega$
is fixed. Note also the following "quenched" theorem in
\cite{FuSh99}: for any $f$ in $L^2_0$, for any $\epsilon>0$, for
$\PP$-almost every $\omega$,
$$\displaystyle {{1 \over \sqrt{n}}} \sum_{k=1}^n f(\tau_k(\omega)
\ldots \tau_1(\omega)\cdot) = o(\log^{3/2+\epsilon}n).\footnote{This
is true in a more general abstract situation (see \cite{FuSh99})}$$

Our main result here is the following:

\vskip 3mm {\bf Theorem} \it Let ${\cal A}$ be a proximal and
totally irreducible finite set \footnote{ The result is still true
if ${\cal A}$ is proximal, irreducible on $\RR^d$ and the semigroup
generated by ${\cal A}$ coincide with the group generated by ${\cal
A}$.} of matrices $d \times d$ with coefficients in $\ZZ$ and
determinant $\pm 1$. Let $\mu$ be a probability measure with support ${\cal A}$ and $\PP =\mu^{\otimes \NN}$ be the product measure on $\Omega:={\cal A}^{\NN}$. Let $f$ be a centered
H\"older function on $\mathbb{T}^d$ or a centered characteristic
function of a regular set. Then, if $f \not \equiv 0$, for
$\PP$-almost every $\omega$ the limit $\sigma(f):= \lim_n {1\over
\sqrt{n}} \Vert S_{n}(\omega,f)\Vert_{2}$ exists and is positive,
and
$$\displaystyle \left ({{1 \over \sigma(f)\sqrt{n}}}
\sum_{k=1}^n f(\tau_k(\omega)\ldots \tau_1(\omega)\cdot) \right)_{n
\geq 1}$$ converges in distribution to the normal law
$\mathcal{N}(0,1)$ with a rate of convergence. \rm

\vskip 3mm An analogous result has been proved for positive $2
\times 2$ matrices and differentiable functions $f$ in
\cite{AyLiSt07} via a different method.

\vskip 3mm The paper is organized as follows. In Section
\ref{freq-tlc} we give sufficient conditions that ensure the
approximation by a normal law of the distribution of the normalized
sums ${1 \over \|S_n\|_2} S_N f$. The proof is based on the method
of multiplicative systems (cf. Koml\`os \cite{Ko73}) (see B. Petit
\cite{Pe92} for an other application of this method). In Section
\ref{exemp-alea}, we address the case of a product of independent
matrices and prove a "quenched" CLT as mentioned above. The key
inequalities are deduced from results of Guivarc'h and Raugi
(\cite{GuRa85}, \cite{Gu90}). Section \ref{sect-statio} is devoted
to the general stationary case under stronger assumptions on the set
$\cal A$, in particular for $2 \times 2$ positive matrices.

\vskip 3mm {\bf Acknowledgements} \  We thank Guy Cohen for his
valuable comments on a preliminary version of this paper.


\tableofcontents

\section{Preliminaries}\label{freq-tlc}


\subsection{A criterion of Koml\`os \label{methodekomlos}}
\vskip 4mm
\goodbreak {\bf Notations} Let $d$ be an integer $\geq 2$ and $\Vert\cdot\Vert$ be the norm
on $\mathbb{R}^{d}$ defined by $\Vert x\Vert=\max_{1\leq i\leq
d}\vert x_{i}\vert, x\in\mathbb{R}^d$. We denote by $d(x,y):=
\inf_{n \in \ZZ^d}(\|x-y -n\|)$ the distance on the torus. The
characters of the torus, $t \rightarrow e^{2\pi i \langle n,
t\rangle}$ for $n=(n_1,...,n_d) \in \ZZ^d$, are denoted by
$\chi(n,t)$ or by $\chi_n(t)$.

The Fourier coefficients of a function $h \in L^2(\mathbb{T}^d)$ are
denoted by $(\hat h(n), n \in \ZZ^d)$. The H\"older norm of order
$\alpha$ of an $\alpha$-H\"older function $f$ on the torus is
$$\|f\|_{\alpha}=\|f\|_{\infty}+\sup_{x\neq
y}{{|f(x)-f(y)|}\over{d(x,y)^{\alpha}}}.$$

A subset $E$ of the torus is said to be regular if there exist $C>0$
and $\alpha\in]0,1]$ such that
$$\lambda(\{t\in\mathbb{T}^d: \
d(t,\partial E)\leq \varepsilon\})\leq C\varepsilon^\alpha, \forall
\varepsilon > 0.$$

The Lebesgue measure on $\mathbb{T}^d$ is denoted by $\lambda$. It
is invariant under the action of automorphisms of the torus. The
action of a product of automorphisms $\tau_j... \tau_i$, $j \geq i$,
corresponds to the action on the characters of the matrices which
define the automorphisms {\it by composition on the right side}. If
$A_1, A_2, \dots$ is a sequence of matrices, if $i\leq j$ are two
positive integers, we use the notation
\begin{eqnarray}
A_i^j := A_i \dots A_j. \label{Aij}
\end{eqnarray}

\vskip 3mm {\bf Multiplicative systems}

In the proof of the central limit theorem we will use a lemma on
"multiplicative systems" (cf. Koml\`os \cite{Ko73}). The
quantitative formulation of the result will yields a rate of
convergence in the CLT. The proof of the lemma is given in appendix.

\begin{lem}\label{komlos}
Let $u$ be an integer $\geq 1$, $(\zeta_{k})_{0\leq k\leq u-1}$ be a
sequence of length $u$ of real bounded random variables, and $a$ be
a real positive number. Let us denote, for $x \in \RR$:
\begin{eqnarray*}
Z(x) &=& \exp({ix}\sum_{k=0}^{u-1}\zeta_{k}(.)),\
Q(x,.) = \prod_{k=0}^{u-1}(1+{ix}\zeta_{k}(.)),\\
Y &=& \sum_{k=0}^{u-1}\zeta_{k}^{2}, \ \ \delta = \max_{0\leq k \leq
u-1} \|\zeta_{k}\|_\infty.\end{eqnarray*}

There is a constant $C$ such that, if $|x| \, \delta\leq 1$, then
\begin{eqnarray} |\EE[Z(x)]
- e^{- {1\over 2}a \,x^2}| \leq C u |x|^3 \delta^3 + {1\over 2} x^2
\|Q(x)\|_2 \|Y - a \|_2 + |1-\EE [Q(x)]|. \label{majKomlos}
\end{eqnarray}
If moreover $|x| \|Y - a\|_2^{1\over 2} \leq 1$, then
\begin{eqnarray}
|\EE[Z(x)]  - e^{- {1\over 2}a \,x^2}| \leq C\,  u \,|x|^3 \delta^3
+ (3 + 2\,e^{-{1\over 2} a \, x^2} \, \|Q(x)\|_2) |x|
\|Y-a\|_2^{1\over 2} + e^{-{1\over 2} a x^2} |1 - \EE [Q(x)]|. \
\label{majKomlos2}
\end{eqnarray}

If $\EE[Q(x,.)] \equiv 1$, the previous bound reduces to
\begin{equation}
C u \,|x|^3 \delta^3 + (3 + 2 \,\,e^{-{1\over 2} a \, x^2} \,
\|Q(x)\|_2) \, |x| \|Y-a\|_2^{1\over 2}. \label{majKomlos3}
\end{equation}
\end{lem}

\vskip 4mm \goodbreak \subsection{Bounding $|\EE [e^{i x {S_n \over
\|S_n\|_2}}] - e^{-{1\over 2} x^2}|$}

\vskip 5mm The application of the lemma to the action on the torus
of matrices in $SL(d, \mathbb{Z})$ requires a property of
"separation of the frequencies" which is expressed in the following
property.
\begin{ppty} {\it Let $(A_1,\ldots, A_n)$ be a finite set of
matrices in $SL(d,\mathbb{Z})$ and let $(D,\Delta)$ be a pair of
positive reals.  We say that the property ${\cal S}(D,\Delta)$ is
satisfied by $(A_1,\ldots, A_n)$ if the following property is
satisfied:

Let $s$ be an integer $\geq 1$. Let $1\leq \ell_1 \leq \ell_1' \leq
\ell_2 \leq \ell_2' \leq \ell_3\leq ... < \ell_s \leq \ell_s'\leq n$
be any increasing sequence of $2s$ integers, such that $\ell_{j+1}
\geq \ \ell_j' + \Delta$ for $j=1,..., s-1$. Then for every families
$p_1, p_2, ..., p_s$ and $p_1', p_2', ..., p_s'$ $\in \ZZ^d$ such
that $A_1^{\ell_s'} p_s' + A_1^{\ell_s} p_s \not = 0$ and $\|p_j\|,
\|p_j'\| \leq D$ for $j=1,..., s$, we have:

\begin{equation}
\sum_{j=1}^s [A_1^{\ell_j'} p_j'  + A_1^{\ell_j} p_j] \not = 0.
\end{equation}}\end{ppty}


A particular case of the property is the following. Let $s$ be an
integer $\geq 1$. Let $\ell_1 < \ell_2 < ... < \ell_s$ be any
increasing sequence of $s$ integers such that $\ell_{j+1} \geq \
\ell_j + \Delta$ for $j=1, ...,s-1$. Then for every family $p_1,
p_2, ..., p_s \in \ZZ^d$ such that $p_s \not = 0$ and $\|p_j\| \leq
D$ for $j=1,..., s$, we have:
\begin{equation}\sum_{j=1}^s A_1^{\ell_j} p_j \not = 0.\end{equation}
This condition implies a multiplicative property as shown by the
following lemma:

\vskip 3mm
\begin{lem} \label{multi}
Let $(D,\Delta)$ be such that the property ${\cal S}(D,\Delta)$
holds with
respect to the finite sequence of matrices $(A_1,\ldots, A_n)$.
Let $g$ be a trigonometric polynomial such that $\hat g(p)
=0$ for $\|p \| > D$. If $\ell_1 < \ell_2 < ... < \ell_s$ is an
increasing sequence of integers such that $\ell_{j+1} \geq \ \ell_j
+ \Delta$ for $j= 1,..., s-1$, then $$\displaystyle {\int \prod_{j =
1}^s g(\tau_{\ell_j}... \tau_1 t) \ dt = 0}.$$
\end{lem}
Recall that the transformation $\tau_{\ell}$ is associated to the
matrix $A_{\ell}$ as said in the introduction.

In what follows,
relative sizes of $D$ and $\Delta$ will be of importance. The
interesting case is when ${\cal S}(D,\Delta)$ is satisfied with
$\Delta$ small compared to $D$ (say $\Delta$ of order $\ln D$). We
will now focus on the characteristic function
$$\EE[e^{i x S_n}]  = \int_{\TT^d} e^{ix \sum_{\ell=0}^{n -1}
g_n(\tau_\ell ... \tau_1 t)} \, dt$$ for a real trigonometric
polynomial $g_n$, where $S_n$ are the ergodic sums
$$S_n(t) = \sum_{\ell=1}^n g_n(\tau_\ell
... \tau_1 t).$$

We will use the inequality given by the following lemma for large
integers $n$. Typically, if $\|S_n\|_2$ is of order $\sqrt{n}$, it
will be applied when $\Delta_n$ is small compared with
$n^{\beta/2}$.
 \vskip 3mm
\begin{lem}\label{bound}
Let $n$ be an integer. Suppose that there exist $\beta\in]0,1[$,
$D_n>0$, $\Delta_n>0$ such that $\Delta_n<{1\over 2} n^\beta$ and
${\cal S}(D_n,\Delta_n)$ is satisfied with
respect to the finite sequence of matrices $(A_1,\ldots, A_n)$. Then, if $g_n$ is a real
trigonometric polynomial with $\hat g_n(p) =0$ for $\|p \|
> D_n$, there exists a polynomial function $C$ with positive
coefficients such that, for $|x|\| g_n \|_\infty n^{\beta}\leq
\|S_n\|_2$  and $|x|\| g_n \|_\infty^{1/2} n^{1+{3\beta} \over 4}
\leq \|S_n\|_2$:
\begin{eqnarray}
&&|\EE [e^{i x {S_n \over \|S_n\|_2}}] - e^{- {1\over 2} x^2}| \nonumber\\
&\leq&  C(\|g_n\|^{1/2}_\infty) [|x| \|S_n\|_2^{-1}\Delta_n n^{1-\beta \over 2}+ |x|^3
\|S_n\|_2^{-3} n^{1 + 2\beta} + |x| \|S_n\|_2^{-1} n^{1+{3\beta} \over
4}  \nonumber \\ && \ \ + |x|^2
\|S_n\|_2^{-1} n^{1-\beta \over 2} \Delta_n+|x|^2 \|S_n\|_2^{-2}
n^{1-\beta}\Delta_n^2]. \label{majFin0}
\end{eqnarray}
\end{lem} \Proof \ \ The proof of (\ref{majFin0}) is given in several steps.

A) {\it Replacement of $S_n$ by a sum with "gaps"}

In order to apply Lemma \ref{komlos}, we replace the sums $S_n$ by a
sum of blocks separated by an interval of length $\Delta_n$.

Let $\beta\in ]0, 1[$, $D_n$, $\Delta_n$ and $g_n$ as in the
statement of the lemma. We set:
\begin{eqnarray}
&& \ v_n:=\lfloor n^\beta \rfloor,
\ u_n:= \lfloor n/v_n \rfloor\leq 2n^{1-\beta}, \label{notatrou1}\\
&&L_{k,n} := k v_n , \ R_{k,n} := (k+1)v_n -\Delta_n, \label{notatrou2}\\
&&I_{k,n} :=[L_{k,n} , \ R_{k,n}], {\rm \ for \ } 0 \leq k \leq u_n
-1. \label{notatrou3}
\end{eqnarray}

Let $S_n'(t)$ be the sum with "gaps":
\begin{equation} S_n'(t) := \sum_{k=0}^{u_n -1}
T_{k,n}(t),
\end{equation}
where \begin{equation} T_{k,n}(t) := \sum_{L_{k,n} < \ell  \leq
R_{k,n}} g_n(\tau_\ell ... \tau_1 t).\end{equation}

The interval $[1,n]$ is divided into $u_n$ blocks of length $v_n -
\Delta_n$. The number of blocks is almost equal to $n^{1-\beta}$ and
their length almost equal to $n^\beta$. The integers $L_{k,n}$ and
$R_{k,n}$ are respectively the left and right ends of the blocks,
which are separated by intervals of length $\Delta_n$. The array of
r.v.'s $(T_{k,n})$ is a "multiplicative system" in the sense of
Koml\`os.

\vskip 4mm \goodbreak {\it Expression of $|T_{k,n}(t)|^2$}

\begin{eqnarray*}
|T_{k,n}(t)|^2 &=& (\sum_{{\ell'} \in I_{k,n}} \sum_{p' \in \ZZ^d}
\hat g(p') \chi(A_1^{\ell'} p', t )) \, (\sum_{\ell \in I_{k,n}}
\sum_{p \in \ZZ^d} \overline {\hat g(p)} \chi(-A_1^{\ell} p, t)) \\
&=& \sum_{p, p' \in \ZZ^d} \ \sum_{\ell, \ell' \in I_{k,n}} \hat
g(p') \overline {\hat g(p)} \chi(A_1^{\ell'} p' -A_1^{\ell} p, t)
\\&=& \sigma_{k,n}^2 + W_{k,n}(t),
\end{eqnarray*}
with
\begin{eqnarray}
\sigma_{k,n}^2 &:=& \int |T_{k,n}(t)|^2 dt = \sum_{p, p' \in \ZZ^d}
\hat g(p') \overline {\hat g(p)} \ \sum_{\ell, \ell' \in I_{k,n}}
1_{A_1^{\ell'} p' = A_1^{\ell} p}, \label{sigmak}\\ W_{k,n}(t)&:=&
\sum_{p, p' \in \ZZ^d}  \hat g(p') \overline {\hat g(p)} \
\sum_{\ell, \ell' \in I_{k,n} : A_1^{\ell'} p' \not = A_1^{\ell} p}
\chi(A_1^{\ell'} p' -A_1^{\ell} p, t). \nonumber
\end{eqnarray}

\vskip 4mm B) {\it Application of Lemma \ref{komlos}} \ We will now
apply Lemma \ref{komlos} to the array of r.v.'s $(T_{k,n}, 0 \leq k
\leq u_n-1)$.  For a fixed $n$, we use the same notations as in the
lemma: $u=u_n$ and for $k=0,...,u_n-1$
\begin{eqnarray*}
\zeta_k &=& T_{k,n}, \  Y = Y_n= \sum_{k=0}^{u_n-1} |T_{k,n}|^2, \\
a &=&a_n = \EE(Y_n) = \sum_k \sigma_{k,n}^2.
\end{eqnarray*}

With the notation of the lemma, the expression of $Q_n(x,t)$ is
\begin{equation}
Q_{n}(x,t)=\prod_{k=0}^{u_n-1}\left(1+i{x} T_{k,n}(t)\right).
\end{equation}
First let us checked that $\EE[Q_n(x,.)] = 1, \forall x$. The
expansion of the product gives
$$Q_{n}(x,t)=1+\sum_{s=1}^{u_n} \left({i x}\right)^s
\sum_{0\leq k_{1}<\dots <k_{s}\leq u_n-1} \prod_{j=1}^s
T_{k_{j},n}(t).$$

The products $\prod_{j=1}^s T_{k_{j},n}(t)$ are combinations of
expressions of the type: $\chi(\sum_{j=1}^s A_1^{\ell_j} p_j, t)$,
with $\ell_j \in I_{k_j, n}$ and $\|p_j\| \leq D_n$. So, by the
property ${\cal S}(D,\Delta)$, $\sum_{j=1}^s A_1^{\ell_j} p_j \not =
0$, and $\displaystyle \int\prod_{j=1}^{s}T_{k_{j},n}(t)\ dt = 0$
(cf. Lemma \ref{multi}).

\vskip 3mm Now we successively bound the quantities involved in
Inequality (\ref{majKomlos3}).

\vskip 3mm B1) {\it Bounding $u_n \delta_n^3$}

$$u_n \delta_n^3 = u_n \max_{0\leq k\leq u_n-1} \|T_{k,n} \|_\infty^3 \leq
C n^{1-\beta} \|g_n\|_\infty^3n^{3\beta} = C \|g_n\|_\infty^3n^{1+ 2 \beta}.$$

B2) {\it Bounding $\|Y_n - a_n\|_2$}

If $U_1, ..., U_L$ are real square integrable r.v.'s such that
$$\EE [(U_k - \EE\,U_k) (U_{k'} - \EE\,U_{k'})] =0, \forall 1 \leq k
<k'\leq L,$$ then the following inequality holds
\begin{eqnarray*}
\| \sum_k U_k - \sum_k \EE\,[U_k]\|_2^2 &=& \sum_k \EE \, [U_k^2] -
(\sum_k \EE \, [U_k])^2\\
&\leq&  \sum_k \EE \, [U_k^2] \leq L \max_k \|U_k\|_\infty \max_k
\EE(|U_k|).
\end{eqnarray*}

We apply this bound to $U_k = (T_{k,n})^2$ and $L = u_n$ (remark
that $T_{k,n}^2 = \sigma_{k,n}^2 + W_{k,n}$ and that we have
orthogonality: $\int W_{k,n} W_{k',n} dt =0$, $ 1 \leq k < k' <u_n$,
due to the choice of the gap. Using rough bounds for $\| T_{k,n}
\|_\infty^2$ and $\| T_{k,n} \|_2^2 $, it implies
$$\|\sum_k T_{k,n}^2 - \sum_k \sigma_{k,n}^2\|_2^2
\leq u_n\| g_n \|_\infty^2    v_n^2 \max_k \sigma_{k,n}^2 \leq 2\|
g_n \|_\infty^2
 n^{1 + \beta} \max_k \sigma_{k,n}^2.$$

So we have
\begin{eqnarray}
|x| \|Y_n - a_n\|_2^{1\over 2} \leq 2^{1/4}\| g_n \|_\infty^{1/2}
|x| n^{1+\beta \over 4} \max_k \sigma_{k,n}^{1\over 2}. \label{xYn}
\end{eqnarray}

B3) {\it Bounding $\EE |Q_n(x)|^2$}

This is the main point. We have
\begin{eqnarray}
|Q_n(x,t)|^2 &=& \prod_{k=0}^{u_n-1} (1+ {x^2} |T_{k,n}(t)|^2) =
\prod_{k=0}^{u_n-1} [1+ {x^2} \sigma_{k,n}^2 + {x^2} W_{k,n}(t)] \\
&=& \prod_{k=0}^{u_n-1} [1+ {x^2} \sigma_{k,n}^2] \
\prod_{k=0}^{u_n-1} [1+ {x^2\over 1 + x^2\sigma_{k,n}^2} W_{k,n}(t)]
\label{Qnx}
\end{eqnarray}

We will show that the integral of the second factor with respect to
$t$ is equal to 1. The first factor in (\ref{Qnx}) is constant and
the bound $1 + y \leq e^y, \forall  y \geq 0$, implies
$$\prod_{k=0}^{u_n-1} [1+ {x^2} \sigma_{k,n}^2]
\leq e^{x^2\sum_{k=0}^{u_n-1} \sigma_{k,n}^2} = e^{a_n x^2}.$$ Hence
the bound
$$\int |Q_n(x,t)|^2 \ dt \leq e^{a_n x^2}.$$

It remains to show that $$\int \ \prod_{k=0}^{u_n-1} [1+ {x^2\over
1+ x^2 \sigma_{k,n}^2} W_{k,n}(t)] \ dt = 1.$$

In the integral the products $W_{k_1} (t)...W_{k_s}(t)$, $0\leq k_1
< ... < k_s < u_n$, are linear combinations of expressions of the
form
$$\chi(\sum_{j=1}^s [A_1^{\ell_j'} p_j' - A_1^{\ell_j} p_j] ,t),$$
where $\ell_j, \ell_j' \in I_{k_j,n}$, $A_1^{\ell_j'} p_j' \not =
A_1^{\ell_j} p_j$, $j=1,...,s$ and $p_j, p_j'$ are vectors with
integral coordinates and norm $\leq D_n$ which correspond to the non
null terms of the trigonometric polynomial $g_n$.

\vskip 3mm As ${\cal S}(D_n,\Delta_n)$ is satisfied, our choice of
gap in the definition of the intervals $I_{k_j,n}$ implies
$\sum_{j=1}^s (A_1^{\ell_j'} p_j' - A_1^{\ell_j} p_j) \not = 0$ and
so the integral of the second factor in (\ref{Qnx}) reduces to 1.

\vskip 3mm From the previous inequalities (in particular (\ref{xYn})
and (\ref{majKomlos3}) of Lemma \ref{komlos}) we deduce that, if
$|x|\| g_n \|_\infty n^{\beta}\leq 1$  and $|x|2^{1/4}\| g_n
\|_\infty^{1/2}
 n^{1+\beta \over 4} \max_k \sigma_{k,n}^{1\over 2}\leq 1$, then
\begin{eqnarray}
|\EE [e^{i {x} S'_n}] - e^{- {1\over 2}a_n x^2}|
& \leq& |x|^3 u_n
\delta_n^3 + (3 + 2\,e^{-{1\over 2} a_n \, x^2} \, \|Q(x)\|_2)\, |x|
\|Y-a_n\|_2^{1\over 2} \nonumber\\
& \leq& C(|x|^3 \|g_n\|_\infty^3n^{1 + 2\beta} + |x|
\|g_n\|_\infty^{1/2}n^{1+\beta \over 4} \max_k \sigma_{k,n}^{1\over
2} ).\label{majtrou}
\end{eqnarray}

\vskip 2mm C) {\it Bounding the difference between $S_n$ and $S'_n$}

\vskip 3mm Recall that $S_n$ is the sum $\sum_{1}^n
g_n(\tau_k...\tau_1 x)$ and $S'_n = \sum_k T_{k,n}$ is the sum with
gaps. We still have to bound the error made when replacing $S_n$ by
$S_n'$:
\begin{eqnarray*} \|S_n - S_n'\|_2^2 &=& \int
|\sum_{k=0}^{u_n -1} \sum_{R_{k,n} <
\ell \leq L_{k+1,n}} g_n(\tau_\ell ... \tau_1 t)|^2 \ dt\\
&=& \sum_{k=0}^{u_n -1} \int |\sum_{R_{k,n} < \ell \leq L_{k+1,n}}
g_n(\tau_\ell ... \tau_1 t)|^2 \ dt \\
&&+ 2 \sum_{0 <k < k' \leq u_n -1} \int \sum_{R_{k,n} < \ell \leq
L_{k+1,n}} g_n(\tau_\ell ... \tau_1 t) \sum_{R_{k',n} < \ell' \leq
L_{k'+1,n}} g_n(\tau_{\ell'} ... \tau_1 t) \ dt.
\end{eqnarray*}

The length of the intervals for the sums in the integrals is
$\Delta_n$. The second sum in the previous expression is 0 by Lemma
\ref{multi} (since $n^\beta-\Delta_n>\Delta_n$). Each integral in
the first sum is bounded by $\|g_n\|_\infty^2 \Delta_n^2$. It
implies:
\begin{equation}
\|S_n - S_n'\|_2^2 \leq \|g_n\|^2_2  \, \Delta_n^2 \, u_n \leq 2
\|g_n\|_\infty^2 n^{1-\beta}\Delta_n^2. \label{erreur-trous}
\end{equation}

Thus, we have
\begin{eqnarray} |\|S_n\|_2^2 - \|S_n'\|_2^2| &\leq& 2
\|S_n\|_2 \|S_n - S_n'\|_2+ \|S_n - S_n'\|_2^2 \nonumber \\ &\leq&
2\sqrt{2} \|S_n\|_2 \|g_n\|_\infty n^{1-\beta \over
2}\Delta_n+2\|g_n\|_\infty^2 n^{1-\beta}\Delta_n^2. \label{majSn2}
\end{eqnarray}

On an other hand, setting $Z_n(x) = e^{ix {S_n}}$, $Z_n'(x) = e^{ix
{S_n'}}$, we have:
\begin{eqnarray}
|\EE [Z_n(x) - Z_n'(x)]| &&\leq \EE [|1 - e^{ix (S_n -S_n')}|] \leq
|x| \EE [|S_n -S_n'|] \leq |x| \|S_n -S_n'\|_2 \nonumber \\
&&\leq \sqrt{2} |x| \|g_n\|_\infty n^{1-\beta \over 2}\Delta_n.
\label{majZn}
\end{eqnarray}

\vskip 4mm D) {\it Conclusion}

Now we gather the previous bounds. Recall that $\sigma^2_{k,n}$
(defined by (\ref{sigmak})) is bounded by
$C\|g_n\|_\infty^2n^{2\beta}$.

\vskip 3mm From (\ref{majZn}), (\ref{majSn2}) and (\ref{majtrou}),
we get that, if $|x|\| g_n \|_\infty n^{\beta}\leq 1$  and $|x|\|
g_n \|_\infty^{1/2} n^{1+{3\beta} \over 4} \leq 1$:
\begin{eqnarray*}
&&|\EE [e^{i {x} S_n}] - e^{- {1\over 2}\|S_n\|_2^2 x^2}| \\
&\leq& |\EE [e^{i {x} S_n}] - \EE [e^{i {x} S_n'}]| + |\EE [e^{i {x}
S_n'}] - e^{- {1\over 2}a_n x^2}| + |e^{- {1\over 2}a_n x^2} -
e^{-{1\over 2}\|S_n\|_2^2 x^2}|\\
&\leq& |\EE [e^{i {x} S_n}] - \EE [e^{i {x} S_n'}]| + C(|x|^3 \|g_n\|_\infty^3n^{1 + 2\beta} +
|x| \|g_n\|_\infty^{1/2}n^{1+\beta \over 4} \max_k \sigma_{k,n}^{1\over 2} )
+ {1\over 2} x^2 |a_n -\|S_n\|_2^2|\\
&\leq&  C(\|g_n\|_\infty) [|x| \Delta_n n^{1-\beta \over 2}+ |x|^3
n^{1 + 2\beta} + |x| n^{1+{5\beta} \over 4}  + |x|^2 \|S_n\|_2
n^{1-\beta \over 2} \Delta_n+ |x|^2 n^{1-\beta}\Delta_n^2].
\end{eqnarray*}
Replacing $x$ by $x  \|S_n\|_2^{-1}$, we obtain Inequality
(\ref{majFin0}) of the lemma. \fdm


\vskip 4mm
\section{Products of independent matrices in $SL(d, \ZZ)$ \label{exemp-alea}}

\subsection{Products of matrices (reminders)}

Let ${\cal A}$ be a finite set of matrices $d\times d$ with
coefficients in $\ZZ$ and determinant $\pm 1$. Let $H$ be the
semi-group generated by ${\cal A}$.

We assume that there is a contracting sequence in $H$ (proximality).
This property holds if ${\cal A}$ contains a matrix with a simple
dominant eigenvalue. We assume also total irreducibility of $H$. It
means that, for every $r$, the action of $H$ on the exterior product
of $\bigwedge_r \RR ^d$ has no invariant finite unions of non
trivial sub-spaces (cf. \cite{Rau97} for this notion).

Let $\mu$ be a probability on ${\cal A}$ such that $\mu(\{A\}) > 0$
for every $A \in {\cal A}$ and let
$$\Omega:={\cal A}^{\NN}= \{\omega = (\omega_n), \omega_n \in {\cal
A}, \forall n \in \NN \}$$ be the product space endowed with the
product measure $\PP=\mu^{\bigotimes \NN}$. For every element
$\omega$ in $\Omega$, we denote by $A_k(\omega)$ (or simply $A_k$)
its $k$-th coordinate. In other words, we consider a sequence of
i.i.d. random variables $(A_k)$ with values in $SL(d, \ZZ)$ and
distribution $\mu$, where $\mu$ is a discrete probability measure
with support ${\cal A}$.

In this section, we prove a central limit theorem with a (small)
rate of convergence for the action of the product $A_n(\omega)\ldots
A_1(\omega)$ on the torus for a.a. $\omega$. This establishes, in a
more general setting and by a different method, a "quenched" central
limit theorem obtained in \cite{AyLiSt07}. The proof relies on
results on products of random matrices obtained by Guivarc'h, Le
Page and Raugi (\cite{GuRa85}, \cite{Gu90}, \cite{LeP80},
 \cite{GuLeP04}). Let us describe the results we will
need.

The group $G=SL(d, \ZZ)$ acts on the projective space ${\bf
P}^{d-1}$. We denote by $(g,x) \rightarrow g.x$ the action. For
$\mu$ a probability measure on $G$, this define a $\mu$-random walk
on ${\bf P}^{d-1}$, where the probability for going from $x$ to
$g.x$ is $d\mu(g)$ (in our case $\mu$ has a finite support  ${\cal
A}$ and $d\mu(g)$ is just $\mu(g)$). Let $X$ and $A$ be two
independent random variables, respectively with values in ${\bf
P}^{d-1}$ and $G$, and with distribution $\nu$ and $\mu$. Then the
distribution of $A.X$ is $\mu\ast\nu$, where
$$(\mu\ast\nu)(\varphi)=\int_{{\bf P}^{d-1}}\int_G \varphi(g.x)\ d\mu(g)\ d\nu(x).$$

The measure $\nu$ is $\mu$-stationary if $\mu\ast\nu=\nu$, i.e.,
when$X$ and $A.X$ have the same distribution. If ${\cal A}$ is
proximal and totally irreducible, then there is a unique
$\mu$-stationary measure on ${\bf P}^{d-1}$ denoted by $\nu$.

\vskip 3mm {\bf Notations} Let $\varphi$ be a function on ${\bf
P}^{d-1}$. We set
$$[\varphi]=\sup_{u,v\in {\bf P}^{d-1},u\neq
v}{{|\varphi(u)-\varphi(v)|}\over{e(u,v)}},$$ where $e$ is the
distance on ${\bf P}^{d-1}$ given by the sinus of the angle between
two vectors. A function $\varphi$ is said to be Lipschitz if
$\|\varphi\|=\|\varphi\|_\infty+[\varphi]<\infty$.

\vskip 3mm For a matrix $M \in Gl(d, \RR)$, we denote by
\begin{eqnarray}
M=N(M)\ D(M)\ K(M) \label{polardec}
\end{eqnarray} its Iwasawa decomposition with $N(M)$ an upper triangular matrix
and $A(M)$ a diagonal matrix with positive diagonal entries.
For $M = A_1(\omega)\ldots A_n(\omega)$, we write
\begin{eqnarray*}&&N^{(n)}(\omega):= N(A_1(\omega)\ldots A_n(\omega)),\\
&&D^{(n)}(\omega):= D(A_1(\omega)\ldots A_n(\omega)),\\
&&K^{(n)}(\omega):= K(A_1(\omega) \ldots A_n(\omega)),
\end{eqnarray*}
so that by (\ref{polardec}): $A_1^n(\omega) := A_1(\omega) \ldots
A_n(\omega) = N^{(n)}(\omega) \, D^{(n)}(\omega) \,
K^{(n)}(\omega).$

\vskip 3mm Let $a_i^{(n)}(\omega):=D_{ii}(A_1^n(\omega))$, $i= 1,
..., d$,  be the diagonal coefficients of the diagonal matrix
$D^{(n)}$.
\begin{prop}
If ${\cal A}$ is proximal and totally irreducible, then there exist $\delta>0$,
$C>0$, and $\rho\in]0,1[$ such that, for every Lipschitz function $\varphi$ on  ${\bf P}^{d-1}$,\\
\begin{eqnarray}
& & \|\EE(\varphi(A_n\ldots A_1.x))-\nu(\varphi)\|_\infty\leq C\rho^n\|\varphi\|,\label{Lepage}\\
& & \int_{\Omega}{\left({a_i^{(n)}(\omega)}\over{a_{i+1}^{(n)}(\omega)}
\right)^\delta}\ d\omega\leq C\rho^n, \ \forall i =1, \cdots, d-1,\label{Guiv1}\\
& &\sup_n\int_{\Omega}\|N^{(n)}(\omega)\|^\delta\ d\PP(\omega)<\infty,\label{Guiv2}\\
& &\sup_{x\in {\bf S}^{d-1}}\int_{{\bf
S}^{d-1}}|\langle x,y\rangle |^{-\delta}d\nu(y)<+\infty\label{Guiv3}.
\end{eqnarray}
\end{prop} \proof \ These statements are consequences of important results of Lepage and
Guivarc'h. They can be deduced: (\ref{Lepage}) from \cite{LeP80},
(\ref{Guiv1}), (\ref{Guiv2}), (\ref{Guiv3}) respectively from
Theorems 5, 6, 7' of \cite{Gu90}. \fdm

\vskip 3mm Let us derive some consequences that we will need.

We first remark that almost surely $a_d^{(n)}(\omega)>
a_i^{(n)}(\omega)$, for $i \in [1, d-1[$, for $n$ large enough.

The Markov inequality and (\ref{Guiv1}) show that there are
constants $C>0$, $\zeta>1$, $\xi_0 \in ]0, 1[$, and a set ${\cal
E}_n$ of measure $\leq C\xi_0^n$ such that, if $\omega$ does not
belong to ${\cal E}_n$, then
\begin{eqnarray}
a_d^{(n)}(\omega)>\zeta^n a_i^{(n)}(\omega), \ \forall i =1, \cdots, d-1.\label{gdsquotients}
\end{eqnarray}

Let $e_d$ be the last element of the canonical basis of $\RR^d$. As
${^t N^{(n)}}$ is lower triangular we have
\begin{eqnarray*}
\EE(\varphi(^tA_n\ldots ^tA_1.e_d))&=&\EE(\varphi(^tK^{(n)}{D^{(n)} }{^t N^{(n)}}.e_d))\\
&=&\EE(\varphi(^tK^{(n)}.e_d)).
\end{eqnarray*}
If we consider the set $\{^tA, A \in \cal A\}$ of transposed
matrices, the conditions of proximality and irreducibility are also
satisfied and we have the same results with a $^t\cal A$-stationary
measure $\nu'$. Because of (\ref{Lepage}) there
exists $\beta_0\in (0,1)$ such that
\begin{eqnarray*}
\|\EE(\varphi(^tA_1\ldots ^tA_n.e_d))-\nu'(\varphi)\|_\infty\leq
C\beta_0^n\|\varphi\|.
\end{eqnarray*}
But
$^tA_n\ldots ^tA_1$ and $^tA_1\ldots ^tA_n$ have
the same distribution. So, for $\varphi$ a function on ${\bf
P}^{d-1}$, under the proximality and irreducibility conditions, there
exists $\beta_0\in (0,1)$ such that
\begin{eqnarray}
\|\EE(\varphi(^tK^{(n)}.e_d))-\nu'(\varphi)\|_\infty\leq
C\beta_0^n\|\varphi\|.\label{loideK}
\end{eqnarray}

\vskip 3mm The probability $\nu'$ satisfies the regularity property
(\ref{Guiv3}). We deduce that there exist $C>0$, $\delta>0$ such
that, for $x\in {\bf S}^{d-1}$ and $\varepsilon>0$, we have
\begin{eqnarray}
 \nu'\{y\in {\bf S}^{d-1}: \ |\langle x,y\rangle
|<\varepsilon\}\leq C \varepsilon^{\delta}.\label{tranche}
\end{eqnarray}
We are now ready to begin our proof.


\subsection{Separation of frequencies}

\vskip 3mm
\begin{lem} \label{clef} There exist $n_0\in\NN$, $C>0$, $\alpha\in ]0,1[$ and
$\beta \in]0,1[$ such that, if $\varepsilon_n\geq \beta^n$, then for
$n\geq n_0$ and for every vector $x$:
\begin{eqnarray}
\PP (\|A_1\ldots A_nx\|\leq \varepsilon_n\|A_1\ldots A_n\|\|x\|)\leq
C\varepsilon_n^\alpha. \label{compar-norm}
\end{eqnarray}
 \end{lem}
\Proof \ \
Let $(e_1,\ldots, e_d)$ be the canonical basis of $\RR^d$. Let $x$ be a unit vector. By expanding $x$ in the
orthonormal basis $(^{t}K^{(n)} e_1,\ldots,^{t}K^{(n)}e_d)$, we
get
$$A_1(\omega)\ldots A_n(\omega)x=\sum_i\langle x, ^{t}K^{(n)}
e_i\rangle \, a_i^{(n)}N^{(n)}e_i.$$

Hence:
$$
\|A_1(\omega)\ldots A_n(\omega)(x)\| \leq \ a_d^{(n)}\
\|N^{(n)}\|\sum_i|\langle x, ^{t}K^{(n)} e_i\rangle| \leq \
a_d^{(n)}\ \|N^{(n)}\|\|x\| \leq \ a_d^{(n)}\ \|N^{(n)}\|.$$

In particular the norm $\|A_1(\omega)\ldots A_n(\omega)\|$ is less
than $a_d^{(n)}\ \|N^{(n)}\|$.

On the other hand we have
$$\|A_1(\omega)\ldots A_n(\omega)(x)\| \geq a_d^{(n)}|\langle x,
^{t}K^{(n)} e_d\rangle|-\|N^{(n)}\|\sum_{i=1}^{d-1}a_i^{(n)}.$$

If the conditions
\begin{eqnarray*}
|\langle {x}, \,^{t}K^{(n)}e_d\rangle |\geq 2\varepsilon_n\|N^{(n)}\|, \
a_d(A_1^n)>\zeta^n a_i(A_1^n), \forall i \in [1, d-1[,
\end{eqnarray*}
are both satisfied, we have
\begin{eqnarray*}
\|A_1(\omega)\ldots
A_n(\omega)(x)\| &\geq& 2\varepsilon_n\|N^{(n)}\| a_d(A_1^n)
\ (1- {{\sum_{i=1}^{d-1}a_i^{(n)}}\over{a_d(A_1^n)\varepsilon_n}})\\
&\geq& 2\varepsilon_n\|N^{(n)}\| a_d(A_1^n)\ (1- (d-1) \zeta^{-n}\varepsilon_n^{-1}).
\end{eqnarray*}

Take $\zeta_0\in(1,\zeta)$. If $\varepsilon_n\geq \zeta_0^{-n}$ and
$(1- (d-1) {\zeta}^{-n_0}{\zeta_0}^{n_0})\geq 1/2$, then for $n \geq
n_0$, we have:
$$\|A_1(\omega)\ldots A_n(\omega)(x)\|\geq\varepsilon_n\|N^{(n)}\| a_d(A_1^n)\ \geq
\varepsilon_n \, \|A_1(\omega)\ldots A_n(\omega)\|.$$

Thus we have obtained that, for $n\geq n_0$,
\begin{eqnarray*}
&&\PP( \|A_1\ldots A_nx\|\leq \varepsilon_n\|A_1\ldots A_n\|\|x\|)
\leq \PP(|\langle {x}, \,^{t}K^{(n)}e_d\rangle |\\
&&\ \ \leq 2\varepsilon_n\|N^{(n)}\|) +\PP(a_d(A_1^n)\leq\zeta^n
a_i(A_1^n)).
\end{eqnarray*}
This inequality gives the following one for every $b_n$. We will
later chose $b_n$ related to $\varepsilon_n$.
\begin{eqnarray*}
&&\PP( \|A_1\ldots A_nx\|\leq \varepsilon_n\|A_1\ldots
A_n\|\|x\|) \leq \PP(|\langle {x}, \,^{t}K^{(n)}e_d\rangle | \\
&&\leq 2\varepsilon_nb_n)+\PP(\|N^{(n)}\|\geq b_n) +
\PP(a_d(A_1^n)\leq\zeta^n a_i(A_1^n)).
\end{eqnarray*}
Let $x \in {\bf P}^{d-1}$ and $\varepsilon>0$. By convolution one
can smooth the indicator function of a "strip" in ${\bf P}^{d-1}$.
There exists $(\varphi_\varepsilon^x)$ a Lipschitz function on the
sphere with values between 0 and 1, such that
$\varphi_\varepsilon^x(y) = 1$ on the set $\{y: |\langle x,y\rangle
|<2\varepsilon\}$, = 0 on the set $\{y:|\langle x,y\rangle
|>3\varepsilon\}$, and such that
$\|\varphi_\varepsilon^x\|<C\varepsilon^{-1}$. Using (\ref{loideK})
and (\ref{tranche}), we have:
\begin{eqnarray*}
&&\PP(\{\omega : \  |\langle  x,^{t}K^{(n)} e_d\rangle
|<2\varepsilon\}) \\ &\leq &\int \varphi_\varepsilon^x(^{t}K(A_1(\omega)
\ldots A_n(\omega))e_d)d\PP (\omega)
\leq \int_{{\bf S}^{d-1}}\varphi_\varepsilon^x(v)d\nu'(v)+C\beta_0^n\varepsilon^{-1}\\
&\leq& \nu'\{v : \ |\langle x,v\rangle
|<3\varepsilon\}+C\beta_0^n\varepsilon^{-1} \leq C
(3\varepsilon)^{\delta}+C\beta_0^n\varepsilon^{-1}.
\end{eqnarray*}
By taking $\varepsilon=\varepsilon_nb_n$, it follows
$$\PP(\{\omega : \  |\langle  x,^{t}K^{(n)} e_d\rangle
|<2\varepsilon_nb_n\}) \\ \leq C
(3\varepsilon_nb_n)^{\delta}+C\beta_0^n(\varepsilon_nb_n)^{-1}.$$

On the other hand, we have by (\ref{gdsquotients})
$$\PP(\{\omega : \  a_d(A_1^n(\omega))\leq\zeta^n a_{d-1}(A_1^n(\omega))\})\leq C\xi_0^n,$$
so that
\begin{eqnarray*}
&&\PP (\|A_1\ldots A_nx\|\leq \varepsilon_n\|A_1\ldots A_n\|\|x\|)
\\
&\leq& \PP ( a_d(A_1^n)\leq\zeta^n a_{d-1}(A_1^n)) +\PP ( |\langle
x,^{t}K^{(n)} e_d\rangle |<2\varepsilon_nb_n) +\PP ( \|N^{(n)}\|\geq b_n)\\
&\leq&
C\xi_0^n+C
(3\varepsilon_nb_n)^{\delta}+C\beta_0^n\varepsilon_n^{-1}b_n^{-1}+Cb_n^{-\delta}.
\end{eqnarray*}

We are looking for $\alpha>0$ such that the following inequalities
hold: $b_n^{-\delta}\leq \varepsilon_n^{\alpha}$, $(b_n
\varepsilon_n)^\delta\leq  \varepsilon_n^{\alpha}$, $\xi_0^n\leq
\varepsilon_n^{\alpha}$ and $\beta_0^n\varepsilon_n^{-1}b_n^{-1}\leq
\varepsilon_n^{\alpha}$. Let us take $\alpha=\delta/3$ and
$b_n=\varepsilon_n^{-1/2}$. Then $b_n^{-\delta}\leq
\varepsilon_n^{\alpha}$ and $(b_n  \varepsilon_n)^\delta\leq
\varepsilon_n^{\alpha}$. If $ \varepsilon_n\geq \beta_0^{n\over
{1\over 2}+\alpha}$ and $ \varepsilon_n\geq\xi_0^{n/\alpha}$ then
the two other inequalities are satisfied.

So by taking $\beta > \max(\zeta_0^{-1},\beta_0^{1\over {{1\over
2}+\alpha}}, \xi_0^{1/\alpha})$, we obtain (\ref{compar-norm}).

Remark that the bound is uniform with respect to $x$ in ${\bf
S}^{d-1}$.\fdm

\vskip 5mm
\begin{coro}\label{cor1} Let $0<\gamma_1<\gamma_2<1$. For
a.a. $\omega$, there exists $L_1(\omega) < +\infty$ such that, for $n
\geq L_1(\omega)$, for every vector $p$ with integral coordinates  and
norm less than $e^{n^{\gamma_1}}$:
$$ \|A_1^np\|\geq e^{-n^{\gamma_2}}\|A_1^n\|\|p\|.$$
\end{coro}
\Proof \ \ Let $n \geq n_0$ such that $e^{-n^{\gamma_2}} \geq
\beta^n$. By Lemma \ref{clef}, for each vector $p$ we have
$$ \PP(\|A_1^np\|\leq
e^{-n^{\gamma_2}}\|A_1^n\|\|p\|) \leq C e^{-\alpha n^{\gamma_2}}.$$
Therefore the probability that there is an integral vector $p$ with
a norm less than $e^{n^{\gamma_1}}$ such that
$$ \|A_1^np\|\leq e^{-n^{\gamma_2}}\|A_1^n\|\|p\|$$
is less then $C e^{dn^{\gamma_1}}e^{-\alpha n^{\gamma_2}}$. This is
the general term of a summable series. We conclude by the
Borel-Cantelli lemma. \fdm

\vskip 3mm
\begin{coro}\label{cor2} For every $M>0$, there exists
$F>0$ such that for a.a. $\omega$, there is $L_2(\omega) < +\infty$
such that, for $n \geq L_2(\omega)$, for every vector $p$ with
integral coordinates  and norm less than $n^M$:
$$ \|A_1^np\|\geq {{1}\over{n^F}}\|A_1^n\|\|p\|.$$
\end{coro}
\Proof \ \ Let $n \geq n_0$ such that $n^{-F}\geq \beta^n$. Thus by
Lemma \ref{clef}, for each vector $p$, for $n\geq n_0$, we have:
$$ \PP(\|A_1^np\|\leq
{{1}\over{n^F}}\|A_1^n\|\|p\|) \leq C {{1}\over{n^{\alpha F}}}.$$
Therefore the probability that there is an integral vector $p$ with
a norm less than $n^M$ such that
$$ \|A_1^np\|\leq {{1}\over{n^F}}\|A_1^n\|\|p\|$$
is less than $C n^{dM}{{1}\over{n^{\alpha F}}}$. If $\alpha F>dM+1$
this is the general term of a summable series. We conclude by the
Borel-Cantelli lemma. \fdm

\begin{lem}\label{clef'} Let $\alpha, \beta \in]0,1[$ given by Lemma
\ref{clef}. There exist $n_0\in\NN$, $C>0$, such that for every unit
vector $x$, every sequence $(\varepsilon_n)$ such that
$\varepsilon_n\geq \beta^n$, every $n$ and every integer $r =
n_0,\cdots, n-n_0$, we have
$$\PP(\|A_1\ldots A_nx\|\leq \varepsilon_r\varepsilon_{n-r}\|A_1\ldots
A_r\|\|A_{r+1}\ldots A_n\|)\leq
C(\varepsilon_r^\alpha+\varepsilon_{n-r}^\alpha).$$
\end{lem}
\Proof \ \ If $\|A_1\ldots A_nx\|>\varepsilon_r\|A_1\ldots
A_r\|\|A_{r+1}\ldots A_nx\|$ and  $\|A_{r+1}\ldots A_nx\|>
\varepsilon_{n-r}\|A_{r+1}\ldots A_n\|$ then $\|A_1\ldots A_nx\|>
\varepsilon_r\varepsilon_{n-r}\|A_1\ldots A_r\|\|A_{r+1}\ldots
A_n\|$. Thus
\begin{eqnarray*}
& & \PP(\|A_1\ldots A_nx\|\leq \varepsilon_r\varepsilon_{n-r}\|A_1\ldots
A_r\|\|A_{r+1}\ldots A_n\|) \\
& &\ \ \ \ \ \ \ \ \ \ \ \ \ \ \leq\PP(\|A_1^nx\|\leq \varepsilon_r\|A_{1}^r\|\|A_{r+1}^nx\|)+
\PP(\|A_{r+1}^nx\|\leq \varepsilon_{n-r}\|A_{r+1}^n\|)
\end{eqnarray*}
Lemma \ref{clef} shows that if $\varepsilon_{n-r}\geq \beta^{n-r}$,
then, for $n-r\geq n_0$, we have
$$
\PP(\|A_{r+1}^nx\|\leq \varepsilon_{n-r}\|A_{r+1}^n\|)\leq
C\varepsilon_{n-r}^\alpha.
$$
The bound obtained in Lemma \ref{clef} is uniform in $x$ and the
matrices $A_{r+1}^n$ and $A_{1}^r$ are independent. Thus, if
$\varepsilon_{r}\geq \beta^{r}$ and $r\geq n_0$, one
has
\begin{eqnarray*}
\PP(\|A_{1}^nx\|\leq
\varepsilon_{r}\|A_1^r\|\|A_{r+1}^nx\|)&=&\PP(\|A_{1}^rA_{r+1}^nx\|\leq
\varepsilon_{r}\|A_{1}^r\|\|A_{r+1}^nx\|)\\
& \leq& \int_\Omega \PP(\|A_{1}^r y\|\leq
\varepsilon_{r}\|A_{1}^r\|\|y\|)\ d\PP_{A_{r+1}^nx}(y) \\ & \leq&
C\varepsilon_{r}^\alpha.
\end{eqnarray*}
\fdm

\vskip 4mm
\begin{coro}{\label{corcle}}
Let $\kappa\in(0,1)$ and $F > {2 \over \kappa \alpha}$. For a.e.
$\omega$, there is $L_3(\omega) < +\infty$ such that, for $n \geq
L_3(\omega)$, for every integer $r$ between $n^\kappa$ and
$n-n^\kappa$: \begin{eqnarray} \|A_1\ldots A_n\|\geq
{{1}\over{r^F(n-r)^F}} \|A_1\ldots A_r\|\|A_{r+1}\ldots A_n\|.
\label{surMult} \end{eqnarray}
\end{coro}
\Proof \ \ If $n$ is large enough and $n^\kappa \leq r \leq n -
n^\kappa$, the inequalities $r\geq n_0$ and $r^{-F}\geq\beta^r$ are
satisfied. Let $x$ be a unit vector. By Lemma \ref{clef'}, for every
integer $r$ between $n^\kappa$ and $n-n^\kappa$, we have
$$ \PP(\|A_1\ldots A_nx\|\leq r^{-F}(n-r)^{-F}\|A_1\ldots
A_r\| \, \|A_{r+1}\ldots A_n\|)\leq C(r^{-\alpha F}+(n-r)^{-\alpha
F}).$$

Therefore the probability that there is an integer $r$ between
$n^\kappa$ and $n-n^\kappa$ such that
$$\|A_1\ldots A_nx\|\leq r^{-F}(n-r)^{-F}\|A_1\ldots
A_r\|\|A_{r+1}\ldots A_n\|$$ is less than $2Cn(n^{-\kappa\alpha
F})$. By the choice of $F$ we have $\sum_{n=1}^\infty
nn^{-\kappa\alpha F} < \infty.$

From the Borel-Cantelli lemma, we deduce that, for a.a. $\omega$,
there is $L_3(\omega) < +\infty$ such that, for every $n\geq
L_3(\omega)$, every integer $r$ between $n^\kappa$ and $n-n^\kappa$:
$$ \|A_1\ldots A_nx\|\geq {{1}\over{r^F(n-r)^F}}
\|A_1\ldots A_r\|\|A_{r+1}\ldots A_n\|.$$

In particular this imples (\ref{surMult}). \fdm

\begin{lem}\label{Allr}
There exists $C_1>0$ such that, for a.a. $\omega$, there is
$L_4(\omega) < +\infty$ such that, for $n \geq L_4(\omega)$, for every
integer $\ell \in [1, n]$, for every $r \in [ C_1 \log n, n]$,
$$\|A_\ell^{\ell+r}\|\geq \zeta^{r(d-1)/d}.$$
\end{lem}

\Proof \ \ We have $\|A_1^r\|\geq |a_d^{(r)}|.$ By a previous
result, (cf. (\ref{gdsquotients}), there exist $\zeta>1$ and
$\xi_0\in]0,1[$ such that
$$\PP(\{|a_d^{(r)}|<\zeta^r |a_i^{(r)}| \}) \leq C\xi_0^r, \forall i =1,
\cdots, d-1.$$

On the other hand the product of the $d_i^{(r)}$'s is one. We thus
have:
\begin{equation}
\PP(\|A_1^r\| < \zeta^{((d-1)/d)r})\leq C\xi_0^r.
\end{equation}
As the probability measure $\PP$ is invariant by the shift, for
every integer $\ell$ we have
\begin{equation}
\PP(\|A_\ell^{\ell+r}\|< \zeta^{((d-1)/d)r})\leq C\xi_0^r.
\end{equation}
The probability that there exist $1\leq \ell \leq n$, $C_1\log n\leq
r\leq n$ such that $\|A_\ell^{\ell+r}\|< \zeta^{((d-1)/d)r}$ is
bounded by $Cn^2\xi_0^{C_1\log n}$. If $C_1>-3/\log \xi_0$ the
sequence $Cn^2\xi_0^{C_1\log n}$ is summable. We conclude by the
Borel-Cantelli lemma. \fdm

\vskip 3mm The next proposition on separation of frequencies shows
that, for $M > 1$ and $\gamma \in ]0,1[$, for a.e. $\omega$, for $n$
big enough, the property ${\cal S}(D_n,\Delta_n)$ is satisfied with
respect to the finite sequence of matrices $(A_1(\omega),\ldots,
A_n(\omega))$ for $D_n = n^M$, $\Delta_n = n^\gamma$. It will enable
us to use Inequality (\ref{majFin0}) of Lemma \ref{bound} with a
well chosen $\gamma$.

\begin{prop} \label{sepa2} For every $\gamma\in ]0,1[$, every $M \geq 1$,
for a.e. $\omega$, there exists a rank $L_5(\omega)$ such that, for
every $n\geq L_5(\omega)$, the property ${\cal S}(n^M, n^\gamma)$ is
satisfied with respect the finite sequence of matrices
$(A_1(\omega),\ldots, A_n(\omega))$. That is:  Let $s$ be an integer
$\geq 1$. Let $1\leq \ell_1 \leq \ell_1' \leq \ell_2 \leq \ell_2'
\leq ... \leq \ell_s \leq \ell_s'\leq n$ be an increasing sequence
of $2s$ integers such that $\ell_{j+1} \geq \ \ell_j' +n^\gamma$,
for $j=1,..., s-1$. Then, for every $p_1, p_2, ..., p_s$ and $p_1',
p_2', ..., p_s'$ $\in \ZZ^d$ such that $A_1^{\ell_s'} p_s'  +
A_1^{\ell_s} p_s \not = 0$ and $\|p_j\|, \|p_j'\| \leq n^{M}$ for
$j=1,..., s$, we have:
\begin{equation}
\sum_{j=1}^s [A_1^{\ell_j'} p_j'  + A_1^{\ell_j} p_j] \not = 0.
\label{equal0}
\end{equation}
\end{prop} \Proof \ \ We will use Corollary \ref{corcle} and the gap between the
$\ell_j$'s to obtain a contradiction from the equality
\begin{equation}A_1^{\ell_s'} p_s'  + A_1^{\ell_s}
p_s = -\sum_{j=1}^{s-1} [A_1^{\ell_j'} p_j'  + A_1^{\ell_j} p_j].
\label{equalFreq}\end{equation}

Let us consider two cases:

\vskip 3mm  1) Assume that $\ell_s' - \ell_s$ is small: $0 \leq
\ell_s' - \ell_s \leq n^\eta$ (for some $\eta\in(0,\gamma^2)$) or
that $p'_s=0$.

Write $q_s = A_{\ell_s+1}^{\ell_s'} p_s' + p_s$. It is a non zero
element of $\ZZ^d$ by assumption and its norm is bounded by
$n^{M}\max\{\|A\|\ : \ A\in {\cal A}\}^{n^\eta}\times n^M\leq
2n^MR^{n^\eta}$. We have $\ell_s> n^\gamma$. So the norm of $q_s$ is
bounded by $2\ell_s^{M/\gamma}R^{\ell_s^{\eta/\gamma}}$. According
to Corollary \ref{cor1}, for every $\eta'>\eta/\gamma$, almost
surely, for $n$ large enough, we have
\begin{equation} \|A_1^{\ell_s}q_s\|\geq e^{-\ell_s^{\eta'}}\|A_1^{\ell_s}\|\geq e^{-n^{\eta'}}\|A_1^{\ell_s}\|.
\label{majAq}
\end{equation}
According to Corollary \ref{corcle}, we have:
\begin{eqnarray*}\label{majo}
\|\sum_{j=1}^{s-1} [A_1^{\ell_j'} p_j'  + A_1^{\ell_j} p_j]\|&\leq&
n^M \sum_{j=1}^{s-1} [\|A_1^{\ell'_j}\|+\|A_1^{\ell_j}\|]\\
&\leq& n^M   \left[\sum_{j=1}^{s-1}
{{\|A_1^{\ell_s}\|}\over{\|A_{\ell'_{j+1}}^{\ell_s}\|}}\ell_j^{'F}
(\ell_s-\ell'_j)^F+\sum_{j=1}^{s-1}
{{\|A_1^{\ell_s}\|}\over{\|A_{\ell_{j+1}}^{\ell_s}\|}}
\ell_j^F(\ell_s-l_j)^F\right]
 \\
&\leq& n^{M+2F} \|A_1^{\ell_s}\| \left[\sum_{j=1}^{s-1}
{{1}\over{\|A_{\ell'_{j+1}}^{\ell_s}\|}}+\sum_{j=1}^{s-1}
{{1}\over{\|A_{\ell_{j+1}}^{\ell_s}\|}}\right]
 \\
&\leq& n^{M+2F} \|A_1^{\ell_s}\|  \left[\sum_{j=1}^{s-1}
\zeta^{(d-1)(\ell'_j-\ell_s)/d}+\sum_{j=1}^{s-1}
\zeta^{(d-1)(\ell_j-\ell_s)/d}\right].
\end{eqnarray*}
The last inequality holds because of Lemma \ref{Allr}. As
$\ell_j-\ell_s\geq (j-s)n^\gamma$ and $\ell'_j-\ell_s\geq
(j-s)n^\gamma$, it implies:
\begin{equation}
\|\sum_{j=1}^{s-1} [A_1^{\ell_j'} p_j'  + A_1^{\ell_j} p_j]\|\leq
Cn^{M+2F} \zeta^{-(d-1)n^\gamma/d} \|A_1^{\ell_s}\|.\label{minAl}
\end{equation}

If we take $\gamma$, $\eta$ and $\eta'$ such that $\eta/\gamma
<\eta' <\gamma$, the inequalities (\ref{majAq}) and (\ref{minAl})
show that (\ref{equalFreq}) is not satisfied for large $n$.

\vskip 3mm 2) Now assume that $\ell_s' - \ell_s$ is large: $ \ell_s'
- \ell_s \geq n^\eta$ and $p'_s\neq0$.

On the one hand, we have (Corollary \ref{cor2})
$$\|A_1^{\ell_s'} p_s'\|\geq {{1}\over{n^F}}\|A_1^{\ell_s'}\|.$$
On the other hand, by Corollary \ref{corcle} and Lemma \ref{Allr}:
\begin{eqnarray}
\|A_1^{\ell_s} p_s \|\leq \|p_s\|
{{\|A_1^{\ell'_s}\|}\over{\|A_{\ell_{s+1}}^{\ell'_s}\|}}\ell_s^F(\ell'_s-\ell_s)^F\leq
n^{M+2F} \zeta^{-(d-1)n^\eta/d} \|A_1^{\ell'_s}\|, \label{majz}
\end{eqnarray} thus by (\ref{minAl}) and (\ref{majz}):
$$
\|A_1^{\ell_s} p_s + \sum_{j=1}^{s-1} [A_1^{\ell_j'} p_j'  +
A_1^{\ell_j} p_j]\|\leq Cn^{M+2F}
(\zeta^{-(d-1)n^\gamma/d}+\zeta^{-(d-1)n^\eta/d})
\|A_1^{\ell'_s}\|.$$

In this case as above, Equality (\ref{equalFreq}) does not hold for
$n$ large enough. \fdm

\vskip 3mm The proposition implies that, for every increasing
sequence of $s$ integers, $\ell_1 < \ell_2 < ... < \ell_s<n$ with
$\ell_{j+1} \geq \ \ell_j + n^\gamma$, for $j=1, ...,s-1$, every
$p_1, p_2, ..., p_s \in \ZZ^d$ such that $p_s \not = 0$ and $\|p_j\|
\leq n^{M}$ for $j=1,..., s$, we have:
\begin{equation}\sum_{j=1}^s A_1^{\ell_j} p_j \not = 0;\end{equation}

\begin{lem}\label{dilat2}
There exist $\zeta_1>1$, $\xi_1\in]0,1[$ and $C>0$ such that
\begin{equation}
\PP(\{\forall p\in \ZZ^d,  \|p\|\leq \zeta_1^n\ : \
\|A_1^np\|>\zeta_1^n\})\geq 1-C\xi_1^n.\label{dilat'}
\end{equation}
\end{lem}
\Proof \ \ We have $\|A_1^n\|\geq |a_d^{(n)}|.$ As we have seen
(Lemma \ref{Allr}), there exist $\zeta>1$ and $\xi_0\in]0,1[$ such
that
\begin{equation}
\PP(\|A_1^n\|\geq \zeta^{((d-1)/d)n})\geq 1-C\xi_0^n. \label{dilat1}
\end{equation}
According to Lemma \ref{clef}, if $\xi_2$ is in $]\beta,1[$, there
exist $C>0$ and $\xi_3\in]0,1[$ such that:
\begin{equation}
\PP(\|A_1^np\|\geq \xi_2^n\|A_1^n\|\|p\|)\geq 1-C\xi_3^n.
\label{dilat1'}
\end{equation}
From (\ref{dilat1}) and (\ref{dilat1'}) we deduce that, if $p$ is an
integral vector, there exist $C>0$, $\zeta_2>1$ and $\xi_4\in]0,1[$
such that (take $\xi_2^{-1}< \zeta^{((d-1)/d)}$):
\begin{equation}
\PP(\|A_1^np\|\geq \zeta_2^n)\geq\PP(\|A_1^np\|\geq
\zeta_2^n\|p\|)\geq 1-C\xi_4^n,
\end{equation}
or equivalently $\PP(\|A_1^np\|\leq \zeta_2^n)\leq C\xi_4^n$.

Let $\zeta_1$ be a real number in $]1,\zeta_2[$ such that
$\zeta_1^d\xi_4<1$. By taking the sum in the previous inequality
over integral vectors $p$ in the ball centered at zero of radius
$\zeta_1^n$, we obtain
$$\PP(\{\exists p\in \ZZ^d, \|p\|\leq \zeta_1^n\ : \  \|A_1^np\|\leq
\zeta_2^n\})\leq C\zeta_1^{dn}\xi_4^n.$$

The lemma follows since $\zeta_2>\zeta_1$. \fdm

The following lemma will be used in the approximation of a function
by a trigonometric polynomial. It can be proved by taking the
sequence $(\varphi_n)$ of the products of the Fej\`er kernels in
each coordinate.

\begin{prop} There exist a constant $C>0$ and a sequence of trigonometric
polynomials $(\varphi_n)$ of order less than $d \, n$, such that,
for every $\alpha$-H\"older function $f$ on the torus,
$$\|\varphi_{n}*f-f\|_{\infty}<C\|f\|_{\alpha}n^{-\alpha},$$
and for every $\alpha$-regular subset $A$ of the torus,
$\|\varphi_{n}*1_A-1_A\|_2<Cn^{-\alpha}$.
\end{prop}


\subsection{Variance and CLT}

We denote by $\theta$ the left shift on $\Omega$:
\begin{eqnarray*}\theta :
(A_k(\omega))_{k\geq 1}\longmapsto(A_{k+1}(\omega))_{k\geq
1},\end{eqnarray*}

$\tau_{A_1}(\omega)$ the map on the torus
\begin{eqnarray*}
\tau_{A_1(\omega)} : x\longmapsto A_1(\omega)x,
\end{eqnarray*}

$\theta_\tau$ the transformation on $\Omega\times\TT^d$:
\begin{eqnarray*}
\theta_\tau : ((A_k(\omega))_{k\geq 1},x)\longmapsto
((A_{k+1}(\omega))_{k\geq 1},\tau_{A_1(\omega)}x),
\end{eqnarray*}
and let
$$S_{n}(\omega,f)(x):=\sum_{k=1}^nf(\tau_{A_k(\omega)}\ldots\tau_{A_1(\omega)}x).$$

\begin{prop}
Let $f$ be a H\"older function on the torus not a.e. null and with
zero mean. Then for $\PP$-almost every $\omega\in\Omega$ the
sequence $({n}^{-{1\over 2}}\Vert S_{n}(\omega,f)\Vert_{2})$ has a
limit $\sigma(f)$ which is positive and does not depend on $\omega$.
\end{prop}
\Proof \ \ Denoting $F(\omega, t) :=f(t)$, we have
$S_{n}(\omega,f)(t) = \sum_{k=0}^{n-1}F(\theta_\tau^k(\omega,t))$
and
\begin{eqnarray}
{1\over n} \Vert S_{n}(\omega,f)\Vert^2= \Vert f\Vert^2+ {2\over
n}\sum_{r=1}^{n-1}\sum_{\ell=0}^{n-1-r} \int_{\mathbb{T}^d}
F(\theta_\tau^\ell(\omega,t))F(\theta_\tau^{\ell+r}(\omega,t))\, dt.
\label{sigma_n-omega}
\end{eqnarray}

Under our hypotheses, the "global variance" exists: the following
convergence holds
$$\lim_n {1\over n} \int \int |S_{n}(\omega,f)|^2 \ dt \,
d\omega = \sigma(f)^2.$$ Actually this holds for every centered
function $f$ in $L^2_0(\TT^d)$, and we have \begin{eqnarray}
\sigma(f)^2 = \sum_{r=1}^{\infty}\int_{\Omega\times\mathbb{T}^d} (F
\ F\circ\theta_\tau^{r})(\omega,t)\, d t\ d\omega \label{serDecor}
\end{eqnarray}
as a consequence of the existence of a spectral gap for the operator
of convolution by $\mu$ on $L^2_0(\TT^d)$, which implies the
convergence of the series. Moreover it can be shown that $\sigma(f)
> 0$ if $f$ is not a.e. null (cf. \cite{Gu06}, see also \cite{FuSh99}).
We have to prove that (\ref{sigma_n-omega}) has the same limit
$\sigma(f)^2$ for a.e. $\omega$.

\vskip 3mm  Let us first consider the sum
$${2\over n}\sum_{r=1}^{n^\alpha}\sum_{\ell=0}^{n-1-r}\int_{\mathbb{T}^d}
F(\theta_\tau^\ell(\omega,t))F(\theta_\tau^{\ell+r}(\omega,t))\, d
t.$$

The second term of the right-hand part of the following equality
\begin{eqnarray*}
&&{2\over n}\sum_{r=1}^{n^\alpha}\sum_{\ell=0}^{n-1-r}
\int_{\mathbb{T}^d} F(\theta_\tau^\ell(\omega,t))
F(\theta_\tau^{\ell+r}(\omega,t))\, d t  = \\
&&2\sum_{r=1}^{n^\alpha}\int_{\mathbb{T}^d} {1\over
n}\sum_{\ell=0}^{n-1} (F \
F\circ\theta_\tau^{r})(\theta_\tau^{\ell}(\omega,t))\, d t
-2\sum_{r=1}^{n^\alpha}{1\over
n}\sum_{\ell=n-r}^{n-1}\int_{\mathbb{T}^d}
F(\theta_\tau^{l}(\omega,t))F(\theta_\tau^{\ell+r}(\omega,t))\, d t\\
\end{eqnarray*}
is bounded by $2\|f\|_2^2n^{2\alpha-1}$. For $\alpha<1/2$, it
suffices to consider the first term.

\vskip 3mm Let us denote by $\psi_j$ the function defined on
$\Omega$ by
\begin{equation}
\psi_j=\int_{\TT^d}F \ F\circ \theta_\tau^j \ dt
-\int_{\Omega\times\TT^d}F \ F\circ \theta_\tau^j \ dt \, d\omega.
\label{psij}
\end{equation}

It only depends on the $j$ first coordinates of $\omega$, so that
$$\int_{\Omega}\psi_j \ \psi_j\circ\theta^l\ \, d\PP(\omega) = 0,
{\rm \ if \ } l>j.$$

We claim that, for every $0<\alpha<1$, $j<n^{\alpha}$,
$\eta\geq2\alpha$,
$$\EE[(\sum_{m=0}^{n-1}\psi_j\circ\theta^m)^4] <Cj^2n^{2}<Cn^{2+\eta}.$$

To prove it, let us expand the fourth power of
$\sum_{m=0}^{n-1}\psi_j\circ\theta^m$ :
$$(\sum_{m=0}^{n-1}\psi_j\circ\theta^m)^4=\sum_{i,k,l,m} \psi_j
\circ\theta^i \, \psi_j \circ\theta^k \, \psi_j\circ \theta^l \,
\psi_j\circ\theta^m.$$

The number of 4-uples $(i,k,l,m)$ such that $i<k\leq i+j$ and
$l<m\leq l+j$ is less than $j^2n^2$ and, if $i+j<k$ or $l+j<m$, then
the integral of the corresponding term is equal to zero. For every
$\varepsilon>0$, the probability
$$\PP(\sup_{j=1,\ldots,n^\alpha}|\sum_{m=0}^{n-1}
\psi_j \circ\theta^m| > n^{\beta}\varepsilon)$$ is less than
$\sum_{j=1}^{n^\alpha}\EE((\sum_{m=0}^{n-1}\psi_j\circ\theta^m)^4)/
\varepsilon^4n^{4\beta}$, therefore it is less than
$Cn^{2+\eta+\alpha-4\beta}$. We can chose $\alpha$, $\beta$, $\eta$
such that this sequence is summable and $\alpha<1-\beta$. Let us
take $\beta=0.8$, $\alpha=0.01$ and $\eta=0.02$).

Then almost surely, we have:
\begin{equation} \lim_{n}{{1}\over{n^{\beta}}}\sup_{j=1,\ldots,n^\alpha}|
\sum_{l=0}^{n-1}\psi_j\circ\theta^l|=0. \label{majBeta}
\end{equation} We deduce from (\ref{psij}) that:
\begin{eqnarray*}
&&\sum_{r=1}^{n^\alpha}\int_{\mathbb{T}^d} {1\over
n}\sum_{\ell=0}^{n-1} (F \ F\circ\theta_\tau^{r})
(\theta_\tau^{\ell}(\omega,t))\, d t \\
&=&\sum_{r=1}^{n^\alpha}{1\over n} \sum_{\ell=0}^{n-1}
\psi_r\circ\theta^l +
\sum_{r=1}^{n^\alpha}\int_{\Omega\times\mathbb{T}^d} {1\over
n}\sum_{\ell=0}^{n-1} (F \ F\circ\theta_\tau^{r}) (\theta_\tau^{\ell} (\omega,t))\, d t\ d\omega\\
&=& \sum_{r=1}^{n^\alpha}n^{\beta-1}{1\over
n^\beta}\sum_{\ell=0}^{n-1}\psi_r\circ\theta^l +
\sum_{r=1}^{n^\alpha}\int_{\Omega\times\mathbb{T}^d} (F \
F\circ\theta_\tau^{r})(\omega,t)\, d t\ d\omega.
\end{eqnarray*}

By (\ref{majBeta}) both sequences
\begin{eqnarray*}&&\sum_{r=1}^{n^\alpha}\int_{\mathbb{T}^d} {1\over
n}\sum_{\ell=0}^{n-1-r} (F \
F\circ\theta_\tau^{r})(\theta_\tau^{\ell}(\omega,t))\, d t, \\
&&\sum_{r=1}^{n^\alpha}\int_{\Omega\times\mathbb{T}^d} (F \
F\circ\theta_\tau^{r})(\omega,t)\, d t\ d\omega \end{eqnarray*}
converge toward the sum of (\ref{serDecor}).

\vskip 3mm We now consider the sum
$$\sum_{r=n^\alpha+1}^{n-1}\int_{\mathbb{T}^d} {1\over
n}\sum_{\ell=0}^{n-1-r} (F \
F\circ\theta_\tau^{r})(\theta_\tau^{\ell}(\omega,t))\, d t.$$

For a given integer $n$, let $g_n$ be a centered polynomial of
degree less than $n^M$ such that $\Vert f-g_n\Vert_{2}<n^{-4}$.

\vskip 3mm Corollary \ref{corcle} shows that, almost surely, the
frequencies of the polynomial $g_n\circ\theta_\tau^{l+r}$ are
greater than $r^{-F}l^{-F}\|A_1^l\|\|A_{l+1}^{l+r}\|$, so greater
than $n^{-2F}\|A_1^l\|\zeta^r$. The norms of the frequencies of the
polynomial $g_n\circ\theta_\tau^{l}$ are less than $n^M\|A_1^{l}\|$.
So, if $r>n^\alpha$, almost surely for a sufficiently large $n$,
$\int_{\mathbb{T}^d} {1\over n}\sum_{\ell=0}^{n-1} (g_n \ g_n
\circ\theta_\tau^{r}) (\theta_\tau^{\ell}(\omega,t))\, d t = 0$.

Thus, we have almost surely in $\omega$
\begin{eqnarray*}
\lim_n {\Vert S_{n}(\omega,f)\Vert^2\over n}&=& \Vert f\Vert^2+
2\sum_{r=1}^{\infty}\int_{\Omega\times\mathbb{T}^d}
(F \ F\circ\theta_\tau^{r})(\omega,t)\, d t\ d\omega \\
&=&\lim_n \int_{\Omega\times\mathbb{T}^d} \, {1\over n} \,
[\sum_{k=0}^{n-1}F(\theta_\tau^k(\omega,t))]^2 \ dt d\omega.
\end{eqnarray*} \fdm

\vskip 3mm Remark that the statement of the previous proposition
holds if $f$ is a centered characteristic function of a regular set
of positive measure of the torus.


\vskip 3mm The proof the CLT for $S_n(\omega,f)$ is now an
application of the method of Section \ref{freq-tlc}.

\begin{thm} Let ${\cal A}$ be a proximal and totally irreducible finite
set of matrices $d \times d$ with coefficients in $\ZZ$ and
determinant $\pm 1$. Let $f$ be a centered H\"older function on
$\TT^d$ or a centered characteristic function of a regular set.
Then, if $f \not \equiv 0$, for almost every $\omega$ the limit $\sigma(f)= \lim_n {1\over \sqrt{n}} \Vert
S_{n}(\omega,f)\Vert_{2}$ exists and is positive, and
$$\displaystyle ({{1 \over \sigma(f)\sqrt{n}}}\sum_{k=1}^nf(\tau_k(\omega)\ldots
\tau_1(\omega)\cdot))_{n \geq 1}$$ converges in distribution to the
normal law $\mathcal{N}(0,1)$ with a rate of convergence.
\end{thm}
\Proof \ \ Recall that the transformation $\tau_k$ is the action on
the torus defined by the transposed  matrice of $A_k$.

Let $f$ be a centered H\"older function or a centered characteristic
function of a regular set such that $\sigma(f) \not = 0$. We have
shown above that almost surely $\|S_nf\|_2$ is equivalent to
$\sigma(f)\sqrt{n}$. It suffices to prove convergence of ${S_nf
\over \|S_nf\|_2}$ towards the normal law $\mathcal{N}(0,1)$.

There exists an integer $M$ such that, for every $n$, there is a
trigonometric polynomial $g_n$ of degree less than $n^M$, such that
$$\|S_nf-S_ng_n\|_2\leq n^{-4}.$$
Therefore $|\EE [e^{i x {S_nf \over \|S_nf\|_2}}] - \EE [e^{i x
{S_ng_n \over \|S_ng_n\|_2}}]|$ tends to 0 and ${\|S_ng_n\|_2 \over
\sqrt{n}}$ tends to $\sigma(f) \not = 0$. Almost surely for $n$ big
enough, the norm $\|S_ng_n\|_2$ is greater than ${1\over 2}
\sigma(f)n^{1/2}$.

\vskip 3mm Now we use the notations introduced in Subsection
\ref{methodekomlos}. According to Proposition \ref{sepa2}, we can
apply Inequality (\ref{majFin0}) of Lemma \ref{bound} with
$\Delta_n=n^{\gamma}$ (this is possible if $\gamma < \beta$): for
$|x|\| g_n \|_\infty n^{\beta}\leq \|S_n\|_2$  and $|x|\| g_n
\|_\infty^{1/2} n^{1+{3\beta} \over 4} \leq \|S_n\|_2$
\begin{eqnarray}
&&|\EE [e^{i x {S_ng_n \over \|S_ng_n\|_2}}] - e^{- {1\over 2} x^2}| \nonumber\\
&\leq&  C(\|g_n\|_\infty) [|x| n^{-{\beta  \over 2}}+ |x|^3 n^{-(-2\beta +1/2)} +
|x| n^{-{(-3\beta +1)\over 4}}
 \ + |x|^2 n^{-{\beta  \over 2}} n^\gamma+ |x|^2
n^{-{\beta }}n^{2\gamma}]. \label{majFin3}
\end{eqnarray}

Here the sequence $(\|g_n\|_\infty)$ is bounded. So, by taking
$\beta$ and $\gamma$ such that $0<2\gamma<\beta<1/4$, we obtain that
$|\EE [e^{i x {S_ng_n \over \|S_ng_n\|_2}}] - e^{- {1\over 2} x^2}|$
tends to 0 for every $x$.

\vskip 3mm Using Esseen's inequality it can be shown that there is
at least a rate of convergence of order $n^{-1/40}$ in the CLT (cf.
\ref{vitesse} in the next section). \fdm

\vskip 5mm
\begin{rems} {\rm  1) If ${\cal A}$ reduces to a single ergodic matrix $A$, the system
is clearly not totally irreducible. Nevertheless, as it is well
known, the CLT holds in this case.

\vskip 3mm 2) Let us take two matrices $A$ and $A^{-1}$ with a non
uniform probability, then the system is not totally irreducible, but
we can show that the quenched CLT holds.

\vskip 3mm 4) If we take two matrices $A$ and $A^{-1}$ with equal
probability 1/2, then the CLT does not hold for the global system.
We get a sort of $T,T^{-1}$ transformation and another limit theorem
(see \cite{LeB06}). This makes us think that CLT for a.a. sequence
of matrices could not be true. } \end{rems}

 \vskip 6mm
\section{Stationary products, matrices in $SL(2, \ZZ^+)$ \label{sect-statio}}

In this section we consider the case of a sequence $(A_k)$ generated
by a stationary process. This is more general than the independent
stationary case, but we have to assume the rather strong Condition
\ref{dilEnt-unif} below, a condition which is satisfied by matrices
in $SL(2, \ZZ^+)$. In this case some information about the
non-nullity of the variance can also be obtained. We will express
the stationarity by using the formalism of skew products.

\vskip 2mm
\goodbreak \subsection{Ergodicity, decorrelation}
\vskip 3mm

We consider an ergodic dynamical system $(\Omega, \mu, \theta)$,
where $\theta$ is an invertible measure preserving transformation on
a probability space $(\Omega, \mu)$. We denote by $X$ the torus
$\mathbb{T}^{d}$ and by $\tau_A$ the automorphism of $X$ associated
to a matrix $A \in SL(d, \ZZ)$. Let ${\cal A}$ be a finite set of
matrices in $SL(d, \ZZ)$.

\begin{notas} {\rm Let $\omega \rightarrow A(\omega)$ be a measurable
map from $\Omega$ to ${\cal A}$, and $\tau$ the map $\omega
\rightarrow \tau(\omega)= \tau_{A(\omega)}$. The skew product
$\theta_\tau$ is defined on the product space $\Omega \times X$
equipped with the product measure $\nu := \mu \times \lambda$ by
$$\theta_\tau :\Omega\times X\rightarrow\Omega\times X;\ (\omega,t)\mapsto
(\theta\omega,\tau(\omega)t).$$

Let $F$ be a function in $L^2(\Omega \times X)$ and, for $p \in
\ZZ^d$, let $F_p(\omega)$ be its Fourier coefficient of order $p$
with respect to the variable $t$. $F$ can be written:
$$F(\omega, t)= \sum_{p \in \ZZ^d} F_p(\omega) \chi(p,t),$$ with
$\sum_{p\in \ZZ^d} \int |F_p(\omega)|^2 \ \, d\mu(\omega)  <
\infty$.

\vskip 3mm Let $\mathcal{H}_{\alpha}^{0}$ be the set of
$\alpha$-H\"older functions on the torus with null integral. This
notation is extended to functions $f(\omega, t)$ on $\Omega \times
X$ which are $\alpha$-H\"older in the variable $t$, uniformly with
respect to $\omega$.}
\end{notas}

For $k\geq 1$, $j \geq i$, $\omega\in\Omega$,
$f\in\mathrm{L}^{2}(X,\mathbb{R})$, we write
\begin{eqnarray*}
\tau(k,\omega)&=&\tau(\theta^{k-1}\omega)\dots \tau(\omega), \\
A_i^j(\omega) &=& A(\theta^i \omega)A(\theta^{i+1} \omega)...
A(\theta^{j} \omega), \\ S_{n}(\omega,f)(t)&=&\sum_{k=1}^{n}f
(\tau(k,\omega)t).
\end{eqnarray*}

\vskip 3mm {\it In what follows in this subsection and in subsection
\ref{nonnul}, we assume the following condition \ref{dilEnt-unif}
which implies an exponential decay of correlation:}

\begin{cond}\label{dilEnt-unif}
There are constants $C>0$, $\delta>0$ and $\lambda>1$ such that
$$\forall r \geq 1, \forall  A_{1},..., A_{r} \in {\cal A},
\, \forall p \in\mathbb{Z}^d \setminus \{0\}, \|A_{1}... A_{r} p \|
\geq {C \Vert p \Vert^{-\delta}}\lambda^{r}.$$
\end{cond}

\begin{prop}\label{ergo} Under Condition \ref{dilEnt-unif},
the system $(\Omega\times X, \theta_\tau, \mu\otimes \lambda)$ is
mixing on the orthogonal of the subspace of functions depending only
on $\omega$. For the skew product map the mixing property holds with
an exponential rate on the space of H\"olderian functions. If
$(\Omega, \mu, \theta)$ is ergodic, then the dynamical system
$(\Omega\times X, \theta_\tau, \mu\otimes \lambda)$ is ergodic.
\end{prop} \Proof \ \ Let $G$ be in $L^2(\Omega\times X)$ a trigonometric
polynomial with respect to $t$ for every $\omega$ and such that
$G(\omega,t) = \sum_{0<\|p\| \leq D} G_p(\omega) \chi(p,t)$, for a
real $D \geq 1$. We have:
\begin{eqnarray*}
\langle G\circ \theta_\tau^n, G \rangle_{\nu} &=& \int\int \ (\sum_p
G_p(\theta^n \omega) \chi(A_1^n(\omega)p,t)) (\overline{\sum_q
G_q(\omega) \chi(q,t)})\ dt \, \, d\mu(\omega)  \\
&=& \sum_{p,q} \int \int \ G_p(\theta^n \omega)
\chi(A_1^n(\omega)p,t)) (\overline{ G_q(\omega) \chi(q,t)})\ dt \, \, d\mu(\omega) \\
&=& \sum_{p,q} \int \, G_p(\theta^n \omega) \,
\overline{G_q(\omega)} \ 1_{A_1^n(\omega) p = q} \, \, d\mu(\omega).
\end{eqnarray*}

According to Condition \ref{dilEnt-unif}, there is a constant $C_1$
not depending on $D$ such that $A_1^n(\omega) p \not = q$, for $n
\geq C_1 \ln D$. Thus we have $\langle G\circ \theta_\tau^n,
G\rangle = 0$, for $n \geq C_1 \ln D$.

\vskip 3mm With a density argument this shows that $\lim_n \langle
G\circ \theta_\tau^n, G \rangle_{\nu} = 0$ for a function $G$ which
is orthogonal in $L^2(\nu)$ to functions depending only on $\omega$
(with an exponential rate of decorrelation for H\"olderian functions
in this subspace). If the system $(\Omega, \mu, \theta)$ is ergodic,
this implies ergodicity of the extension. \fdm

\vskip 2mm We are going to prove that, for a.e. $\omega$, the
sequence $({n}^{-{1\over 2}}\Vert S_{n}(\omega,f)\Vert_{2})$,
converges to a limit. The norm $\Vert S_n(\omega,f)\Vert_2$ is taken
with respect to the variable $t$, $\omega$ being fixed.

\vskip 3mm
\begin{prop}\label{cnsvariancenonnulle} For every $f\in
\mathcal{H}_{\alpha}^{0}(\mathbb{T}^d)$, for $\mu$-a.e.
$\omega\in\Omega$, the sequence $({n}^{-{1\over 2}}\Vert
S_{n}(\omega,f)\Vert_{2})$ has a limit $\sigma(f)$ which does not
depend on $\omega$.

Moreover $\sigma(f)=0$, if and only if $f$ is a coboundary: there
exists $h \in L^{2}(\nu)$ such that
\begin{equation}\label{eqncnsvar}
f(t)=h(\theta\omega,\tau(\omega)t)-h(\omega,t), \ \nu-{\rm a.e.}
\end{equation}\end{prop}
\Proof \ \ The convergence of the sequence of (global) variances
(i.e. for the system $(\Omega\times X,\theta_\tau)$)
$$({n}^{-1} \int
\int |S_{n}(\omega,f)|^2 \ dt \, d\mu(\omega))_{n\geq 1}$$ to an
asymptotic limit variance $\sigma^2$ is a  general property of
dynamical systems, for function with a summable decorrelation. In
this case, we also know that $\sigma= 0$ if and only if $f$ is a
coboundary with a square integrable transfer function.

\vskip 3mm The system $(\Omega\times\mathbb{T}^d,\theta_\tau,\mu
\times dt)$ is ergodic according to Proposition \ref{ergo}.

\vskip 3mm Denoting $F(\omega, t) :=f(t)$, we have
$S_{n}(\omega,f)(t) = \sum_{k=0}^{n-1}F(\theta_\tau^k(\omega,t))$,
hence:
\begin{eqnarray*}
{1\over n} \Vert S_{n}(\omega,f)\Vert_2^2&=& {1\over
n}\sum_{\ell=0}^{n-1}\sum_{\ell'=0}^{n-1}\int_{\mathbb{T}^d}
F(\theta_\tau^\ell(\omega,t))F(\theta_\tau^{l'}(\omega,t))\, d t
\\&=&\Vert f\Vert^2+
{2\over n}\sum_{r=1}^{n-1}\sum_{\ell=0}^{n-1-r}\int_{\mathbb{T}^d}
F(\theta_\tau^\ell(\omega,t))F(\theta_\tau^{\ell+r}(\omega,t))\, d t
\\&=&\Vert f\Vert^2+ 2\sum_{r=1}^{n-1} {1\over n} \int_{\mathbb{T}^d} \sum_{\ell=0}^{n-1}
(F.F\circ\theta_\tau^{r})(\theta_\tau^{\ell}(\omega,t))\, d t\\
&-&2\sum_{r=1}^{n-1} {1\over n} \int_{\mathbb{T}^d}
F(\theta_\tau^{\ell}(\omega,t))
\sum_{\ell=n-r}^{n-1}F(\theta_\tau^{\ell+r}(\omega,t))\, d t.
\end{eqnarray*}

\vskip 3mm Condition \ref{dilEnt-unif} insures, for a constant $C$
and for a real $\kappa < 1$, the following inequality:
\begin{eqnarray}
&&|\int_{\mathbb{T}^d} F(\theta_\tau^{\ell }(\omega,t))
F(\theta_\tau^{\ell+r}(\omega,t))\, d t| \nonumber\\&& =
|\int_{\mathbb{T}^d} f(t) \ f(A(\theta^{\ell+r} \omega) ...
A(\theta^{\ell +1} \omega)t)\, d t|
 \leq C \Vert f\Vert_{2}\Vert f\Vert_{\alpha}
\kappa^r. \label{maj-unilr}
\end{eqnarray}

This implies: \begin{eqnarray*}
\left\vert\sum_{r=1}^{n-1}\int_{\mathbb{T}^d} F(\theta_\tau^{\ell
}(\omega,t)){1\over
n}\sum_{l=n-r}^{n-1}F(\theta_\tau^{\ell+r}(\omega,t))\, d
t\right\vert&\leq& C \Vert f\Vert_{2}\Vert f\Vert_{\alpha} \,
{1\over n}\sum_{r=1}^{n-1}r\kappa^r,
\end{eqnarray*}
hence this term tends to 0 if $n \rightarrow +\infty$ and the
convergence of ${1\over n} \Vert S_{n}(\omega,f)\Vert_2^2$ reduces
to that of \begin{equation}\Vert f\Vert^2+ 2\sum_{r=1}^{n-1} {1\over
n} \int_{\mathbb{T}^d} \sum_{\ell=0}^{n-1}
(F.F\circ\theta_\tau^{r})(\theta_\tau^{\ell}(\omega,t))\, d t.
\label{omega-variance} \end{equation}

\vskip 3mm For $\mu$-a.e. $\omega$, for every $r$, by the ergodic
theorem  \begin{eqnarray*} &&\lim_{n\rightarrow+\infty} {1\over
n}\sum_{\ell=0}^{n-1} \int_{\mathbb{T}^d} F(\theta_\tau^{\ell
}(\omega,t)) F(\theta_\tau^{\ell+r}(\omega,t))\, d t \\ &=& \lim_n
{1\over n}\sum_{\ell=0}^{n-1} \int_{\mathbb{T}^d} f(t) \
f(A(\theta^{\ell+r} \omega) ... A(\theta^{\ell +1} \omega)t)\, d t =
\int_{\Omega\times\mathbb{T}^d}(F.F\circ\theta_\tau^{r})\, d \omega
\, dt.\end{eqnarray*}

According to (\ref{maj-unilr}) we can take the limit for $\mu$-a.e.
$\omega$ in (\ref{omega-variance}):
$$\lim_{n\rightarrow+\infty} {1\over n} \Vert S_{n}(\omega,f)\Vert^2 = \Vert
f\Vert_{2}^{2}+2\sum_{r=1}^{+\infty}
\int_{\Omega\times\mathbb{T}^d}(F.F\circ\theta_\tau^{r})\, d \omega
\, dt = \lim_n \int \int {1\over n} |S_{n}(\omega,f)|^2  \ dt \, \,
d\mu(\omega) .$$ \fdm

\vskip 3mm
\begin{rem} {\rm \label{uniq-erg}
The previous proof shows that for a uniquely ergodic system
$(\Omega, \mu, \theta)$ defined on a compact space $\Omega$ (for
instance an ergodic rotation on a torus), the convergence of the
variance given in Proposition \ref{cnsvariancenonnulle} holds for
every $\omega \in \Omega$, if the map $\tau$ is continuous outside a
set of $\mu$-measure 0.}
\end{rem}

\vskip 3mm
\subsection{Non-nullity of the variance \label{nonnul}}

Now we consider more precisely the condition of coboundary. For $j,p
\in \ZZ^d$, we denote by $D(j,p,\omega)$ the set $\{k \geq 0 :
A_0^k(\omega) j = p \}$ and by $c(j,p,\omega):=\#D(j,p,\omega)$. (By
convention, $A_0^0(\omega) = Id$.) We will use the following simple
lemma:

\vskip 3mm \begin{lem} \label{finiJ} Under Condition
\ref{dilEnt-unif}, $\sup_{j \in J, p \in \ZZ^d} \ c(j,p,\omega) <
\infty$, for every finite subset $J$ of $\mathbb{Z}_*^{d}$.
\end{lem} \Proof \ \  Let $j$ be in $J$ and let $k_1 := \inf\{k \in
D(j,p,\omega)\}$. If $k_2$ belongs to $D(j,p,\omega)$ with $k_2 > k_1$, then
$A_1^{k_2}(\omega) j = p = A_1^{k_1} (\omega) j$, so that:
$A_{k_1+1}^{k_2}(\omega) j = j$. According to Condition
\ref{dilEnt-unif}, this implies that the number of such integers
$k_2$ is finite and bounded independently of $p$. As $J$ is finite,
the result follows. \fdm

\vskip 4mm
\begin{prop} Assume Condition \ref{dilEnt-unif}.
Let $f$ be a trigonometric polynomial in $L^2(\mathbb{T}^d)$. If
there exists $g\in L^2(\Omega\times \mathbb{T}^d)$ such that $\int
g\ d\nu=0$ and $f=g-g\circ\theta_\tau$, then $g$ is also a
trigonometric polynomial.
\end{prop} \Proof  \ \ Let
$\displaystyle f=\sum_{j\in J}f_{j}\chi_{j}$, where $J$ is a finite
subset of $\mathbb{Z}^{d}$. Let $g$ be in  $L^2$ such that $\int g \
d\nu = 0$ and $f(t) = g(\theta \omega, \tau(\omega) t) -
g(\omega,t)$.

\vskip 3mm The coboundary relation implies $\sum_{k=0}^{N-1} (1-{k
\over N}) f\circ \theta_\tau^k = g - {1 \over N} \sum_1^N g\circ
\theta_\tau^k$. As $g$ belongs to $L^2$, by ergodicity we deduce the
convergence in $L^2$-norm
$$g = \lim_N \sum_{k=0}^{N-1} (1-{k \over N}) f\circ
\theta_\tau^k,$$ with \begin{equation}\label{SNf} \sum_{k=0}^{N-1}
(1-{k \over N}) f\circ \theta_\tau^k = \sum_{p \in \ZZ^d} \
\sum_{k=0}^N \ [\sum_{j \,: \,A_0^k(\omega) j = p} (1-{k \over N}) \
f_j] \chi_p,
\end{equation}

 Moreover it is known that
the maximal function $\sup_N {1 \over N} |\sum_1^N g\circ
\theta_\tau^k|$ is square integrable.  Therefore, by Fubini, for
a.e. $\omega$, there is $M(\omega) < \infty$ such that
\begin{equation} \sup_N \sum_{p \in
\ZZ^d} \ |\sum_{k=0}^N \ [\sum_{j \,: \,A_0^k(\omega) j = p} (1-{k
\over N}) \ f_j]|^2 < M(\omega) \label{majSomCar} .\end{equation}

\vskip 3mm If $N$ goes to $\infty$, the expression $\sum_{k=0}^N \
[\sum_{j \,: \,A_0^k(\omega) j = p} (1-{k \over N}) \ f_j]$ tends to
the finite sum $\sum_{j \in J} c(j,p,\omega) \ f_j$ (cf. Lemma
\ref{finiJ}). According to (\ref{majSomCar}), by restricting first
the sums to a finite set of indices $p$ and passing to the limit
with respect to $N$ in $\sum_{p} \ |\sum_{k=0}^N \ [\sum_{j \,:
\,A_0^k(\omega) j = p} (1-{k \over N}) \ f_j]|^2$, we obtain finally
$$\sum_{p \in \ZZ^d} \ |\sum_{j \in J} c(j,p,\omega) \ f_j|^2 < M(\omega).$$

For every $p$, as $J$ is finite and as $c(j,p,\omega)$ takes
integral bounded values according to Lemma \ref{finiJ},
$\displaystyle(\vert\sum_{j\in J} c(j,p, \omega) f_{j}\vert)_{p \in
\ZZ^d}$ take only a finite number of distinct values. Let $V$ be the
set of these values and $\delta>0$ a lower bound of $V \setminus
\lbrace 0 \rbrace$.

\vskip 3mm We have $\delta^2 \ \#\lbrace p \in \mathbb{Z}^d :
\sum_{j\in J}c(j,p,\omega)f_{j}\neq 0 \rbrace \leq M(\omega)$, so
that the cardinal is finite for a.e. $\omega$. This shows that $g$
is a trigonometric polynomial. \fdm

\vskip 3mm
\begin{coro} If $f$ is a coboundary and has non negative Fourier
coefficients, then $f(x) = 0$ a.e.
\end{coro} \Proof \ \ By using the fact that
$c(j,p,\omega) \in \NN$, we get:
\begin{eqnarray*}
\Vert g(\omega, .)\Vert_{2}^{2} &=&\sum_{p}\left(\sum_{j\in
J}c(j,p,\omega)f_{j}\right)^2 \geq \sum_{p}\left(\sum_{j\in
J}c(j,p,\omega)f_{j}^{2}\right)
\\ &\geq &\sum_{j\in J}\left(\sum_{p}c(j,p,\omega)\right) f_{j}^{2}.
\end{eqnarray*}
For $j\neq0$, we have $\displaystyle \sum_{p}c(j,p,\omega)=+\infty$,
which implies $f_{j}=0$. \fdm

\vskip 3mm  The previous results  allow to obtain a "quenched" CLT
(i.e. for a.e. $\omega$) in the stationary case for positives
matrices in $SL(2,\ZZ)$, with (for trigonometric polynomials) a
criterion of non-nullity of the variance. Moreover, when the Fourier
coefficients of $f$ are nonnegative, then the variance is $ > 0$.


\vskip 3mm
\subsection{${\cal A} \subset SL(2, \ZZ^+)$}

We consider in this subsection a finite set ${\cal A}$ of matrices
in $SL(2,\mathbb{Z}^+)$ with positive coefficients. We study the
asymptotical behavior of the products $A_i^j:= A_i... A_j$, where
$A_i, ... ,A_j$, $i \leq j$, is any choice of matrices in ${\cal
A}$.

\vskip 3mm Let $M$ be a $2 \times 2$ matrix with $>0$ coefficients
and having different real eigenvalues $r= r(M), s = s(M)$, $r>s$.

Let $$\tilde M =\pmatrix {r & 0\cr 0 & s,} \ F = \pmatrix { a & b
\cr c & d }$$ be respectively the diagonal matrix conjugate to $M$
and the matrix such that $M = F \tilde M F^{-1}$, with $ad-bc =1$.

\begin{lem} The matrix $M$ can be written:
$$M = \pmatrix {(r-s)u+s & -(r-s)v\cr (r-s)w & -(r-s)u+r},$$
with $u = ad \in ]0, 1[$, $v=ab <0$, $w= cd >0$.
\end{lem} \Proof \ \
The positivity of the coefficients of $M$ implies that $ v <0$,
$w>0$. By multiplying the relation $ad-bc =1$ by $ad$, we obtain
$u^2 -vw =u$, thus $u^2-u=vw<0$. \fdm

\begin{lem} There exist a constant $C$ such that for every $p$ in
$\ZZ_*^2$, and every product $M$ of $n$ matrices taking values in
${\cal A}$, if $n \geq C \ln  \|p\|$, then $Mp \in
\RR_+^2\cup\RR_-^2 $.
\end{lem} \Proof \ \ Let $\lambda := {w\over u} = {u-1 \over v}$.
We have $\lambda > 0$ and we can rewrite the matrix $M$ as
$$M = r\pmatrix {u & \lambda^{-1}(1-u)\cr \lambda u & 1-u}
+ s\pmatrix {1-u & -\lambda^{-1}(1-u)\cr -\lambda u & u}.$$

Thus, for every vector $X = \pmatrix {x\cr y}$, $MX = r (ux  +
\lambda^{-1}(1-u) y)\pmatrix {1\cr \lambda} + s (x- \lambda^{-1}y)
\pmatrix {1-u\cr -\lambda u}$.

The eigenvectors of $M$ are $\pmatrix {1\cr \lambda}$ and $\pmatrix
{1-u\cr -\lambda u}$, corresponding respectively to the eigenvalues
$r$ and $s$.

As $M$ is a product of $n$ matrices of ${\cal
A}$, it maps the cone $\RR_+^2$ strictly into itself:
$$M \RR_+^2  \subset \bigcup_{A \in {\cal A}} A \RR_+^2.$$

It follows that the slope $\lambda = \lambda(M)$ of the positive
eigenvector of $M$ is bounded from below and above by constants
which only depend on ${\cal A}$: there exists $\delta >0$ such that
$\delta \leq \lambda \leq \delta^{-1}$.

Let us write  $r\zeta + \varphi$ and $\lambda r\zeta + \psi$ the
components of $MX$ with:
\begin{eqnarray*}
\zeta := u x + \lambda^{-1} (1-u)y, \  \varphi := s(x - \lambda^{-1}
y) (1-u), \ \psi := - s(x - \lambda^{-1} y) \lambda u.
\end{eqnarray*}

There exist constants $C' >0$ and $\gamma> 1$ such that the positive
eigenvalue $r(M)$, for $M$ a product of $n$ matrices taking values
in ${\cal A}$, satisfies: $r(M) \geq C' \gamma^n$.

\vskip 3mm As $s(M) = r(M)^{-1}$, we  have $s(M) \leq C'^{-1}
\gamma^{-n}$ and, as $\delta \leq \lambda \leq \delta^{-1}$,
$$\max(|\varphi| ,|\psi| )\leq C'^{-1} \delta^{-1}\gamma^{-n}\| X\|.$$

\vskip 3mm Let $X \in \ZZ^2$ be non zero. Up to a replacement of $X$
by $-X$, we can assume that $\zeta \geq 0$. The vector $MX$ having
non zero integer coordinates, we have:
$$r\zeta +|\varphi| + \lambda r\zeta + |\psi| \geq |r\zeta +\varphi|
+ |\lambda r\zeta + \psi| \geq 1.$$

Thus:
$$ r\zeta \geq {1\over 1 + \lambda} - {1\over 1 + \lambda} (|\varphi| + |\psi|),$$
and
$$
r\zeta + \varphi\geq {1\over 1 + \lambda} - {1\over 1 + \lambda}
((2+\lambda)|\varphi| + |\psi|),\ \ \ \ \ \ \lambda r\zeta +
\psi\geq  {\lambda\over 1 + \lambda} - {1\over 1 + \lambda}
(\lambda|\varphi| + (1+2\lambda)|\psi|).
$$
As $\max(|\varphi| ,|\psi| )\leq C'^{-1} \delta^{-1}\gamma^{-n}\|
X\|$, there exists $C>0$ such that if $n\geq C\ln\| p\|$ then
$r\zeta + \varphi > 0$ and $\lambda r\zeta + \psi> 0$ that is $MX\in
\RR_{+*}^2$. \fdm

\vskip 3mm
\begin{coro} Let $(A_k)_{k\geq 1}$ be a sequence of matrices taking values
in ${\cal A}$.  Denote by $\tau_k : x \rightarrow A_k x {\rm \ mod \
} \ZZ^2$ the corresponding automorphisms of the torus. Then, for
almost every $x$ in $\TT^2$, the sequence $(\tau_k... \tau_1 x)_{k
\geq 1}$ is equidistributed in $\TT^2$.
\end{coro}

\vskip 3mm
\begin{coro} \label{dilat} There exist constants
$C_1 >0, \gamma > 1$, and $c$ such that for every $p \in \ZZ^2
\setminus \{0\}$:
\begin{equation}
\|A_1^{\ell+r} p \| \ \geq C_1 \gamma^{r - c\log \|p\|} \|A_1^{\ell}
\|, \ \forall \ell, r \geq 1. \label{dilat-ineq}
\end{equation}
\end{coro}

For vectors $q \in \ZZ_+^2$ belonging to some cone strictly
contained in the positive cone, the norm $\|A^n q\|$ is comparable
to the norm $\|A^n\|$ and there are constants $C> 0$ and $\lambda >
1$ such that $\|A^n q\| \geq C \lambda^n$. Therefore, ${\cal
S}(D,\Delta)$ is a consequence of (\ref{dilat-ineq}).

\begin{coro} \label{sepa} For every $D>0$ there exists $\Delta$
such that ${\cal S}(D,\Delta)$ holds with respect to any products of
matrices in ${\cal A}$.
\end{coro} \Proof \ \  Let us suppose that $\sum_{j=1}^s [A_1^{\ell_j'} p_j'  +
A_1^{\ell_j} p_j] = 0$, i.e.
\begin{equation}
A_1^{\ell_s'} p_s'  + A_1^{\ell_s} p_s = - \sum_{j=1}^{s-1}
[A_1^{\ell_j'} p_j'  + A_1^{\ell_j} p_j]. \label{non-sep}
\end{equation}
Inequality (\ref{dilat-ineq}) ensures inequalities such as:
$$\|A_1^{\ell_j} p_j\| \leq D \|A_1^{\ell_j}\| \leq C_1^{-1} D
\gamma^{-(\ell_s - \ell_j) + c \ln{\Vert q_s\Vert}}
\|A_1^{\ell_s}q_s\|,$$ for $q_s\in\mathbb{Z}^2_*$ and $\Vert
p_j\Vert\leq D$. Then, using the gaps between the $\ell_j$ we will
get a contradiction. More precisely we consider two cases.

\vskip 3mm 1) $\ell_s' - \ell_s$ small: $0 \leq \ell_s' - \ell_s
\leq \rho_1$, where $\rho_1$ will be defined later.

Write $q_s = A_{\ell_s+1}^{\ell_s'} p_s' + p_s$. This is a non-zero
vector in $\ZZ^2$ and its norm is less than $2D\max_{A \in {\cal
A}}\,\|A\|^{\rho_1}$. Let $C_2 := \ln \max_{A \in {\cal A}}\,\|A\|$.
We deduce from (\ref{dilat-ineq}):
\begin{eqnarray*}
\|A_1^{\ell_s'} p_s'  + A_1^{\ell_s} p_s\| &=& \|A_1^{\ell_s} q_s\|
\\ &\leq& C_1^{-1} D [\sum_{j=1}^{s-1} \gamma^{-(\ell_s - \ell_j') + c \ln \|
q_s \|} \|A_1^{\ell_s} q_s \| + \sum_{j=1}^{s-1} \gamma^{-(\ell_s -
\ell_j) + c \ln \| q_s \|}
\|A_1^{\ell_s} q_s \|] \\
&\leq& C_1^{-1} D \gamma^{c \ln \| q_s \|} \|A_1^{\ell_s} q_s \| \
[\sum_{j=1}^{s-1} \gamma^{-(\ell_s - \ell_j')}+ \sum_{j=1}^{s-1}
\gamma^{-(\ell_s - \ell_j)}] \\
&\leq& 4 C_1^{-1}  D' \gamma^{c C_2 \rho_1}\ [\sum_{j=1}^{s-1}
\gamma^{-j \Delta)}] \ \|A_1^{\ell_s} q_s \|\\
&\leq& {4 \over C_1 (1-\gamma^{-\Delta})} D' \gamma^{c C_2 \rho_1 -
\Delta} \ \|A_1^{\ell_s} q_s \|,
\end{eqnarray*}
with $D' = \gamma^{c ln(2D)} D$.

\vskip 3mm \goodbreak 2) $\ell_s' - \ell_s \geq \rho_1$.

We can assume that $p_s' \not = 0$. Otherwise we would have $p_s
\not =0$ and we would consider $\|A_1^{\ell_s} p_s\|$ instead of
$\|A_1^{\ell_s'} p_s'\|$. Still using (\ref{dilat-ineq}) we get:
\begin{eqnarray*}
\|A_1^{\ell_s'} p_s'\| &\leq&  \|A_1^{\ell_s} p_s\| + C_1^{-1} D[
\sum_{j=1}^{s-1} \gamma^{-(\ell_s' - \ell_j') + c \ln \| p_s'\|}
\|A_1^{\ell_s'} p_s' \| + \sum_{j=1}^{s-1} \gamma^{-(\ell_s' -
\ell_j) + c \ln \| p_s'\|}
\|A_1^{\ell_s'} p_s' \|]\\
&\leq& C_1^{-1} D \gamma^{c \ln \| p_s'\|} \|A_1^{\ell_s'} p_s' \|
\, [ \gamma^{-(\ell_{s}' - \ell_s)} + \sum_{j=1}^{s-1}
\gamma^{-(\ell_{s}' - \ell_{j}')} + \sum_{j=1}^{s-1}
\gamma^{-(\ell_s'
- \ell_j)}] \\
&\leq& C_1^{-1} D [\gamma^{c \ln D - \rho_1} + 2{\gamma^{c \ln D -
\Delta} \over  (1-\gamma^{-\Delta})}] \ \|A_1^{\ell_s'} p_s'\|.
\end{eqnarray*}

Chose $\rho_1$ such that $C_1^{-1} D \gamma^{c \ln D - \rho_1} <
{1\over 2}$, then $\Delta$ such that
\begin{eqnarray*}
&&2 C_1^{-1} {D \gamma^{c \ln D - \Delta} \over (1-\gamma^{-\Delta})} < {1\over 2}, \\
&&{4 \over C_1 (1-\gamma^{-\Delta})} D' \gamma^{c C_2 \rho_1 -
\Delta} < 1.
\end{eqnarray*}

The factor in front of  $\|A_1^{\ell_s} q_s \|$ on the right in the
first case is $< 1$ and the factor in front of $\|A_1^{\ell_{s}'}
p_s' \|$ on the right in the second case is $< 1$. In both cases
there is a contradiction. \fdm

\vskip 3mm Corollary \ref{sepa} and Inequality (\ref{majFin0})
enable us to prove a CLT for the action of sequences $A_1^n$. Let
$g_n = g$ be a fixed trigonometric polynomial such that
$\hat{g}(p)=0$ if $\|p\|>D$. Let us take $\Delta$ such that ${\cal
S}(D,\Delta)$ holds (via Corollary \ref{sepa}) and remark that in
the case of $SL(2,\ZZ^+)$ that we are studying the numbers
$\sigma_{k,n}^{1\over 2}$ are bounded by $Cn^{\beta/2}$. Inequality
(\ref{majFin0}) of Lemma \ref{bound} becomes:
\begin{eqnarray}
&&|\EE [e^{i x {S_n \over \|S_n\|_2}}] - e^{- {1\over 2} x^2}| \nonumber\\
&\leq&  C [|x| \|S_n\|_2^{-1} n^{1-\beta \over 2}+ |x|^3
\|S_n\|_2^{-3} n^{1 + 2\beta} + |x| \|S_n\|_2^{-1} n^{1+ 3\beta
\over 4} \nonumber  \\&& \ \ + |x|^2 \|S_n\|_2^{-1} n^{1-\beta \over
2} \Delta+ |x|^2 \|S_n\|_2^{-2} n^{1-\beta}\Delta^2]. \label{majFin}
\end{eqnarray}

If we suppose that $\|S_n\|_2 \geq C n^\delta$, we get:
\begin{eqnarray}
&&|\EE [e^{i x {S_n \over \|S_n\|_2}}] - e^{- {1\over 2} x^2}| \nonumber\\
&\leq&  C [|x| n^{-{(\beta -1 + 2\delta) \over 2}}+ |x|^3
n^{-(-2\beta -1 +3\delta)} + |x| n^{-{(-3\beta -1 +4\delta)\over 4}}
 \nonumber  \\&& \ \  + |x|^2 n^{-{(\beta -1 +2\delta) \over 2}} \Delta+ |x|^2
n^{-({\beta -1 +2\delta)}}\Delta^2]. \label{majFin2}
\end{eqnarray}

{\it Inequality of Esseen}

If $X, Y$ are two r.r.v.'s defined on the same probability space,
their mutual distance in distribution is defined by:
$$d(X,Y) = \sup_{x \in \RR} |\PP(X \leq x) - \PP(Y \leq x)|.$$

Let be $H_{X,Y}(x) := |\EE(e^{i x X}) - \EE(e^{i x Y})|$. Take as
$Y$ a r.v. $Y_\sigma$ with a normal law ${\cal N}(0, \sigma^2)$.

Recall the following inequality (cf. Feller, {\it An introduction to
probability theory and its application}, p. 512): if $X$ has a
vanishing expectation and if the difference of the distributions of
$X$ and $Y$ vanishes at $\pm \infty$, then for every $U>0$,
$$d(X,Y_\sigma) \leq {1\over \pi} \int_{-U}^U H_{X,Y}(x) {dx \over x}
+ {24 \over \pi} {1 \over \sigma \sqrt{2 \pi }} {1\over U}.$$

Taking $X= S_n / \|S_n\|_2$, we have here that $|H_{X,
Y_{\sigma_n}}| \leq \sum_{i=1}^5 n^{-\gamma_i} |x|^{\alpha_i}$,
where the constants are given by (\ref{majFin2}). Thus $d(X,Y_1)$ is
bounded by
$${C \over U} + \sum_{i=1}^5
n^{-\gamma_i} {1\over \alpha_i} U^{\alpha_i}.$$

In order to optimize the choice of $U=U_n$, we take $U_n = n^\gamma$
with $\gamma = \min_i {\gamma_i \over \alpha_i +1}$. This gives the
bound
$$d({S_n \over \|S_n\|_2},Y_1) \leq C n^{-\gamma}.$$

We have to compute
\begin{eqnarray}
\gamma &=& \min({\beta -1 + 2\delta \over 4}, {-2\beta -1 +3\delta
\over 4}, {-3\beta -1 +4\delta \over 8}, {\beta -1 +2\delta \over
6}, {\beta -1 +2\delta \over 3}) \nonumber \\ &=& \min({-2\beta -1
+3\delta \over 4}, {-3\beta -1 +4\delta \over 8}, {\beta -1 +
2\delta \over 6}). \label{beta-delta}
\end{eqnarray}

For $\delta = {1\over 2}$ we get: $\gamma = \min({-4\beta + 1 \over
8}, {-3\beta +1 \over 8}, {\beta \over 6}) =  \min({-4\beta + 1
\over 8},{\beta \over 6})$. Taking $\beta = {3\over 16}$, we obtain
$\gamma = {1\over 32}$. This gives a rate of convergence of order
$n^{-{1\over 32}}$.

\begin{thm}\label{vitesse} Let $(A_k)_{k \geq 1}$ be a sequence
of matrices taking values in a finite set ${\cal A}$ of matrices in
$SL(2, \ZZ_+)$ with $>0$ coefficients. If, for a constant $C_1 > 0$
and a rank $n_0$, $\|S_n\| \geq C_1 n^{1\over 2}$, for $n \geq n_0$,
then for a constant $C$ we have:
\begin{equation} d({S_n \over \|S_n\|_2}, Y_1) \leq C n^{-{1\over
32}}, \forall n \geq n_0. \label{rate1}
\end{equation}
\end{thm}

The previous results can be applied if the limit of $n^{-{1\over
2}}\,\|S_{n}\|_2$ exists and is non zero:  the sequence
$(n^{-{1\over 2}}\,S_{n})_{n \geq 1}$ then tends in distribution
towards the normal law $N(0,1)$ with a rate given by (\ref{rate1}).

\vskip 3mm We can also obtain a rate of convergence of order
$n^{-\delta}$, for some $\delta >0$, for subsequences provided that
the variance $\|S_{n_k}\|_2$ is large enough:

\vskip 3mm Along a subsequence $(n_k)$ such that $\|S_{n_k}\|_2 \geq
C_1 n_k^{\delta}$, with $\delta > 3/7$, the subsequence of
normalized sums $(\|S_{n_k}\|_2^{-1} S_{n_k})$  converges in
distribution towards the normal law ${\cal N}(0,1)$.

\vskip 3mm Indeed, in (\ref{beta-delta}), to obtain a strictly
positive $\gamma$, we have to check the inequalities:
$$-2\beta -1 +3\delta >0, \ -3\beta -1 +4\delta > 0,\  \beta -1 + 2\delta >0.
$$
that is: $$1 - 2\delta < \beta < \min( {3 \delta -1 \over 2},
{4\delta -1\over 3}) = {3 \delta -1 \over 2}.$$

For $\delta > {3\over 7}$ and $\beta = {1\over 7}$, we have $\gamma>
0$.

\vskip 3mm
\begin{rems} {\rm
1) In the previous statements, we have considered the case of
trigonometric polynomials. Using some approximation, it can be
extended to H\"older continuous functions or characteristic functions
of a regular set.

\vskip 3mm 2) If the sequence $(A_n)$ is generated by a dynamical
system $(\Omega,\theta,\mu)$, we have shown that in the case of
$SL(2,\ZZ^+)$-matrices, that either for $\mu$-almost
$\omega\in\Omega$, $(\Vert S_{n}(\omega,f)\Vert_{2})$ is bounded or,
$\mu$-almost $\omega\in\Omega$, the sequence $({n}^{-{1\over
2}}\Vert S_{n}(\omega,f)\Vert_{2})$ has a limit $\sigma(f) > 0$ not
depending on $\omega$. In the later case, the CLT holds.

\vskip 3mm For instance (cf. Remark \ref{uniq-erg}), if the sequence
$(A_n)$ is generated by an ergodic rotation on the circle, with
$A(\omega) = A$ on an interval and $=B$ on the complementary, then
we obtain the CLT for every such sequence.

\vskip 3mm 3) If the dynamical system $(\Omega,\theta,\mu)$ is
weakly mixing, then the characteristic function of a regular set is
never a coboundary for the extended system. Thus we necessarily have
$\sigma(f) > 0$. That is to say that, if $(\Omega,\theta,\mu)$ is
weakly mixing, the CLT holds almost surely for centered
characteristic functions of regular sets.}
\end{rems}


\section{Appendix}

{\bf Proof of Lemma \ref{komlos}}  \ 1) Setting $\psi(y) =
(1+iy)e^{- {1\over2} y^2} e^{-iy}$ and writing $\psi(y) = \rho(y)
e^{i\theta(y)}$, where $\rho(y) = |\psi(y)|$, we have
\begin{eqnarray*}
\ln \rho(y) = {1\over 2}[\ln(1+y^2) - y^2] \leq 0, \ \tan(\theta(y))
= {y - \tan y \over 1 + y \tan y}.
\end{eqnarray*}
An elementary computation gives the following upper bounds for some
constant $C_1$:
\begin{eqnarray}
|\ln \rho(y)| \leq {1\over 4} |y|^4, \  |\theta(y)| \leq C_1 |y|^3,\
\forall y \in [-1, 1]. \label{majo-arg}
\end{eqnarray}

Let us write: $\displaystyle{Z(x) = Q(x) \ \exp(-{1\over 2} x^2 \ Y)
\ [\prod_{k=0}^{u-1} \psi({x\zeta_{k}})]^{-1}}$. Using the fact that
$\ln \rho({x\zeta_{k}})\leq 0$, we have:
\begin{eqnarray*}
|Z(x) - Q(x) \ \exp(-{1\over 2} x^2 \ Y) | &=& |Z(x) - Z(x) \,
\prod_{k=0}^{u-1} \psi({x\zeta_{k}})|
= |1 - \prod_{k=0}^{u-1} \psi({x\zeta_{k}})|\\
&\leq& \displaystyle {|1 - e^{\sum_{k=0}^{u-1} \ln
\rho({x\zeta_{k}})}| + |1 - e^{i\sum_{k=0}^{u-1}
\theta({x\zeta_{k}})}}|\\
&\leq& \displaystyle {\sum_{k=0}^{u-1} |\ln \rho({x\zeta_{k}})| +
\sum_{k=0}^{u-1} |\theta({x\zeta_{k}})}|.
\end{eqnarray*}

If $|x| \delta \leq 1$, where $\delta = \max_k \|\zeta_k\|_\infty$,
we can apply the bound (\ref{majo-arg}). Using the inequality
\begin{equation} |1
- e^s| \leq (e-1) |s| \leq 2 |s|, \forall s \in [-1, 1]
\label{maj-exp} \end{equation} we obtain for a constant $C$:
\begin{eqnarray*}
|Z(x) - Q(x) \ \exp(-{1\over 2} x^2 \ Y) | \leq C {\vert x\vert^{3}
\, \sum_{k=0}^{u-1}\vert \zeta_{k}\vert^3} \leq C u |x|^3 \delta^3.
\end{eqnarray*}

2) Since $Y$ is a positive random variable, we have also:
$$\left\vert \exp(-{1\over 2}x^2 Y)-\exp(- {1\over 2} a \,x^2)\right\vert
\leq {x^2\over2}\left\vert {Y} - a\right\vert.$$

If $|x| \delta \leq 1$ we get:
\begin{eqnarray*}
&&|Z(x) - \exp(-a \,{x^2 \over 2}) Q(x)| \\
&\leq& |Z(x) - Q(x) \, \exp(- {x^2 \over 2}Y)| + |Q(x) \, [\exp(-
{x^2 \over 2}Y) - \exp(-a \,{x^2 \over 2})]|
\\ &\leq& C u \,|x|^3 \delta^3 +{x^2 \over 2} |Q(x)|\, |Y -
a | ;
\end{eqnarray*}
hence, under the condition $|x| \delta \leq 1$, we obtain the upper
bound (\ref{majKomlos}):
\begin{eqnarray*} |\EE[Z(x)]  -
\exp(-a \,{x^2 \over 2})| &=& |\EE[Z(x) - e^{-{1\over 2} a \, x^2}
Q(x) + e^{-{1\over 2} a \, x^2} (Q(x)-1)]|\\
&\leq&|\EE[Z(x) - e^{-{1\over 2} a \, x^2} Q(x)]|
+ e^{-{1\over 2} a \, x^2}|\EE[Q(x)-1)]|\\
&\leq& C u \,|x|^3 \delta^3 +{x^2 \over 2} \EE [|Q(x)|
|Y - a |] + e^{-{1\over 2} a \, x^2} |1-\EE[Q(x)]| \\
&\leq& C u \,|x|^3 \delta^3 +{x^2 \over 2} \|Q(x)\|_2 \|Y - a \|_2 +
|1-\EE[Q(x)]|.
\end{eqnarray*}

3) The bound that we obtain is large in general, because the
integral of $Q(x)$ is of order $e^{{1\over 2}\sigma^2 x^2}$ and the
bound for $Q(x)$ is very large if $x$ is big. If $e^{-{1\over 2} a
\, x^2}\, \|Q(x)\|$ is bounded, we can obtain a more accurate upper
bound.

For $0 \leq \varepsilon \leq 1$, let $A_\varepsilon(x) = \{\omega :
x^2|Y(\omega) - a| \leq \varepsilon \}$. We have the following
bounds:
\begin{eqnarray*}
&&|\EE[1_{A_\varepsilon(x)} (Z(x)  - \exp(-a \,{x^2 \over
2}))]| \\
&\leq& C\,  u \,|x|^3 \delta^3 + \EE[1_{A_\varepsilon(x)}\, (|Q(x)
\, [\exp(- {x^2 \over 2}Y) - e^{-{1\over 2} a \, x^2}]|)]\\ &+&
e^{-{1\over 2} a \, x^2} [|1 - \EE(Q(x)| +
\EE(1_{A_\varepsilon^c(x)} |1 - Q(x)|)]\\
&\leq& C\,  u \,|x|^3 \delta^3 +  e^{-{1\over 2} a \, x^2}
\|Q(x)\|_2
\|1_{A_\varepsilon(x)}[\exp(- {x^2 \over 2}(Y - a )) - 1]\|_2\\
&+& e^{-{1\over 2} a \, x^2} [|1 - \EE(Q(x))| +
\EE(1_{A_\varepsilon^c(x)} |1 - Q(x)|)].
\end{eqnarray*}

From (\ref{maj-exp}) we have
$$\|1_{A_\varepsilon(x)} \, [\exp(- {x^2 \over 2}(Y - a )) - 1]\|_2
\leq 2 \varepsilon,$$ and using Cauchy-Schwarz inequality, we get
$$\EE(1_{A_\varepsilon^c(x)} |1 - Q(x)|)] \leq
\PP(A_\varepsilon^c(x)) + \|Q(x)\|_2
(\PP(A_\varepsilon^c(x)))^{1\over 2},
$$
which implies:
\begin{eqnarray*}
&&|\EE[1_{A_\varepsilon(x)} (Z(x)  - \exp(-a \,{x^2 \over
2}))]| \\
&\leq& C\,  u \,|x|^3 \delta^3 + e^{-{1\over 2} a \, x^2}\|Q(x)\|_2
\, \varepsilon \\&+& e^{-{1\over 2} a \, x^2} [|1 - \EE(Q(x))| +
\PP(A_\varepsilon^c(x)) + \, \|Q(x)\|_2
(\PP(A_\varepsilon^c(x)))^{1\over 2}]
\\ &\leq& C\,  u \,|x|^3 \delta^3 +
e^{-{1\over 2} a \, x^2} \|Q(x)\|_2 \, [\varepsilon +
(\PP(A_\varepsilon^c(x)))^{1\over 2}] + e^{-{1\over 2} a \, x^2} [|1
- \EE(Q(x))| + \PP(A_\varepsilon^c(x))].
\end{eqnarray*}

Choosing $\varepsilon = |x| \|Y - a\|_2^{1\over 2}$, we get
\begin{eqnarray*}
\PP(A_\varepsilon^c(x)) \leq \varepsilon^{-2} x^4 \| Y - a\|_2^2
\leq x^2 \|Y-a\|_2.
\end{eqnarray*}

This yields:
\begin{eqnarray*}
&&|\EE[1_{A_\varepsilon(x)} (Z(x)  - \exp(-a \,{x^2 \over 2}) ]| \\
&\leq& C\,  u \,|x|^3 \delta^3 + 2 |x|\,e^{-{1\over 2} a \, x^2} \,
\|Q(x)\|_2 \, \|Y-a\|_2^{1\over 2} + e^{-{1\over 2} a \, x^2}\, [|1
- \EE(Q(x))| + x^2 \|Y-a\|_2].
\end{eqnarray*}

Thus, assuming $|x| \delta \leq 1 {\rm \ and \ }\  |x| \|Y -
a\|_2^{1\over 2} \leq 1$, we obtain (\ref{majKomlos2}):
\begin{eqnarray*}
&&|\EE[Z(x)]  - \exp(- {1\over 2}a \,{x^2})| \leq
|\EE[1_{A_\varepsilon(x)} (Z(x) - \exp(-a \,{x^2 \over 2}))]| + 2
\PP(A_\varepsilon^c(x) )\\
&\leq& C\,  u \,|x|^3 \delta^3 + 2 \, |x|\,e^{-{1\over 2} a \, x^2}
\, \|Q(x)\|_2 \, \|Y-a\|_2^{1\over 2} + e^{-{1\over 2} a \, x^2}\,
[|1 - \EE(Q(x))|] + 3 \,|x|^2 \|Y-a\|_2.
\end{eqnarray*}
\fdm


Jean-Pierre Conze, Stéphane Le Borgne,  Mikael Roger

conze@univ-rennes1.fr \hfill \break
stephane.leborgne@univ-rennes1.fr \hfill \break
m.mikael.roger@orange.fr

\end{document}